\documentclass[aos,nameyear,dvips]{arximspdf}
\usepackage{graphicx}

\doi{10.1214/10-AOS830}
\volume{39}
\issue{1}
\pubyear{2011}
\firstpage{131}
\lastpage{173}

\makeatletter
\newtheorem{theorem}[definition]{Theorem}
\newtheorem{lemma}[definition]{Lemma}
\newproclaim{remark}[definition]{Remark}
\newproclaim{example}{Example}
\newproclaim{definition}{Definition}
\newproclaim{The non-Markovian-$Y^0$ case}{The non-Markovian-$Y^0$ case}
\newproclaim{The leading-indicator-in-$L$ case}{The leading-indicator-in-$L$ case}
\newtheorem{proposition}[definition]{Proposition}
\makeatother

\begin{document}
\begin{frontmatter}

\title{Causal inference for continuous-time processes when covariates
  are observed only at~discrete~times}
\runtitle{Causal inference for continuous-time processes}

\begin{aug}
\author[A]{\fnms{Mingyuan} \snm{Zhang}\ead[label=e1]{zhangmi@wharton.upenn.edu}},
\author[B]{\fnms{Marshall M.} \snm{Joffe}\thanksref{a1}\ead[label=e2]{mjoffe@mail.med.upenn.edu}}
\and
\author[A]{\fnms{Dylan S.} \snm{Small}\thanksref{a2}\corref{}\ead[label=e3]{dsmall@wharton.upenn.edu}}
\runauthor{M. Zhang, M. M. Joffe and D. S. Small}
\affiliation{University of Pennsylvania}
\thankstext{a1}{Supported in part by NIH Grant CA095415.}
\thankstext{a2}{Supported in part by NSF Grant SES-0961971.}
\address[A]{M. Zhang\\
 D. S. Small \\
 Department of Statistics \\
 The Wharton School\\
  University of Pennsylvania\\
  Philadelphia, Pennsylvania 19104\\
  USA\\
  \printead{e1} \\
  \phantom{E-mail: }\printead*{e3}} 
\address[B]{M. M. Joffe \\
 Department of Biostatistics \\
 \quad\& Epidemiology\\
 University of Pennsylvania \\
 Philadelphia, Pennsylvania 19104\\
 USA \\
 \printead{e2}}
\end{aug}

\received{\smonth{9} \syear{2009}}
\revised{\smonth{4} \syear{2010}}

\begin{abstract}
Most of the work on the structural nested model and
g-estimation for causal inference in longitudinal data assumes a
 discrete-time underlying data generating process.
 However, in some observational studies, it is more reasonable to
assume that the data are generated from a continuous-time
process and are only observable at discrete time points.
When these circumstances arise, the
sequential randomization assumption in the observed discrete-time
data, which is essential in justifying discrete-time g-estimation,
may not be reasonable.
Under a deterministic model, we discuss other useful assumptions that guarantee
the consistency of discrete-time g-estimation.
 In more general cases, when
those assumptions are violated,
we propose a controlling-the-future method that performs at least as well as
g-estimation in most scenarios and which provides consistent estimation in
some cases where g-estimation is severely inconsistent.
We apply the
methods discussed in this paper to simulated data, as well as to a
data set collected following a massive
flood in Bangladesh, estimating the effect of diarrhea on children's
height. Results from different methods are compared in both simulation
and the real application.
\end{abstract}



\end{frontmatter}

\section{Introduction and motivation}\label{s:intro}
In this paper, we study assumptions and methods for making causal
inferences about the effect of a treatment that varies in continuous
time when its time-dependent confounders are observed only at discrete
times. Examples of settings in which this problem arises are given in
Section~\ref{s:intro:example}. In such settings, standard discrete-time
methods such as g-estimation usually do not work, except when certain
conditions are assumed for the continuous-time process. In this paper,
we formulate such conditions. When these conditions do not hold, we
propose a controlling-the-future method which can produce consistent
estimates when g-estimation is consistent and which is still consistent
in some cases when g-estimation is severely inconsistent.

  First, we review the approach of James Robins and collaborators to making
 causal inferences about the effect of a treatment that varies at discrete,
observed times.

\subsection{Review of Robins' causal inference approach for treatments
 varying at discrete, observed times}
\label{s:intro:review}

In a cross-sectional observational study of the effect of a treatment
on an outcome, a usual assumption for making causal inferences is that
there are no unmeasured confounders, that is, that conditional on the
measured confounders, the data is generated as if the treatment were
assigned randomly.  Under this assumption, a consistent estimate of the
average causal effect of the treatment can be obtained from a correct
model of the association between the treatment and the outcome
conditional on the measured confounders [\citet{cochran1965}].  In
a longitudinal study, the analog of the ``no unmeasured confounders''
assumption is that at the time of each treatment assignment, there are
no unmeasured confounders; this is called the \textit{sequential
    randomization} or \textit{sequential ignorability} assumption, given as follows.
 \begin{longlist}
\item[(A1)] The longitudinal data of interest are generated as if the
treatment is randomized in each period, conditional on the current
values of measured covariates and the history of the measured
covariates and the
 treatment.
 \end{longlist}
The sequential randomization assumption implies that decision on
treatment assignment is based on observable history and contemporaneous
covariates, and that people have no ability to see into the future.
Robins (\citeyear{robins1986}) has shown that for a longitudinal study,
unlike for a cross-sectional study, even if the sequential
randomization assumption holds, the standard method of estimating the
causal effect of the treatment by the association between the outcome
and the treatment history conditional on the confounders can provide a
biased and inconsistent estimate.  This bias can occur when we are
interested in estimating the joint effects of all treatment assignments
and when the following conditions hold:
 \begin{longlist}
 \item[(c1)] conditional on past
treatment history, a time-dependent
  variable is a predictor of the subsequent mean of the outcome and
  also a predictor of subsequent treatment;
 \item[(c2)] past treatment history is an independent predictor of the
  time-dependent variable.
  \end{longlist}
Here, ``independent predictor'' means that prior treatment predicts
current levels of the covariate, even after conditioning on other
covariates. An example in which the standard methods are biased is the
estimation of the causal effect of the drug AZT (zidovudine) on CD4
counts in AIDS patients.  Past CD4 count is a time-dependent confounder
for the effect of AZT on future CD4 count since it not only predicts
future CD4 count, but also subsequent initiation of AZT therapy. Also,
AZT history is an independent predictor of subsequent CD4 count [e.g.,
Hern\'{a}n, Brumback and Robins (\citeyear{hernan2002})].

To eliminate the bias of standard methods for estimating the causal
effect of treatment in longitudinal studies where sequential
randomization holds but there are time-dependent confounders satisfying
conditions (c1) and (c2) (e.g., past CD4 counts), Robins
(\citeyear{robins1986}, \citeyear{robins1992}, \citeyear{robins1994},
\citeyear{robins1998}, \citeyear{robins1999}) developed a number of
innovative methods. We focus here on structural nested models (SNMs)
and their associated methods of g-testing and g-estimation.  The basic
idea of the g-test is the following.  Given a hypothesized treatment
effect and a deterministic model of the treatment effect, we can
calculate the \textit{potential outcome} that a subject would have had
if she never received the treatment. Such an outcome is also known as a
\textit{counterfactual outcome}, which is the outcome under a treatment
history that might be contrary to the realized treatment history. If
the hypothesized treatment effect is the true treatment effect, then
this potential outcome will be independent of the actual treatment the
subject received conditional on the confounder and treatment history,
under the sequential randomization assumption (A1). g-estimation
involves finding the treatment effect that makes the g-test statistic
have its expected null value. For simplicity, our exposition focuses on
deterministic rank-preserving structural nested distribution models;
g-estimation also works for nondeterministic structural nested
distribution models.

The SNM and g-estimation were developed for settings in which treatment
decisions are being made at discrete times at which all the confounders
are observed.  In some settings, the treatment is varying in continuous
time, but confounders are only observed at discrete times.

\subsection{Examples of treatments varying in continuous time where covariates are
 observed only at discrete times}\label{s:intro:example}

 \begin{example}[(The effect of diarrhea on children's height)]\label{exam1}
 Diarrheal disease is one of the leading causes of childhood illness in
 developing regions [Kosek, Bern and
 Guerrant (\citeyear{kosek2003})]. Consequently, there is considerable concern
 about the effects of diarrhea on a child's physical and
 cognitive development [Moore et al. (\citeyear{moore2001}), Guerrant et al.
(\citeyear{guerrant2002})].  A data set which provides the opportunity
to study
 the impact of diarrhea on a child's height is a longitudinal
 household survey conducted in Bangladesh in 1998--1999 after
 Bangladesh was struck by its worst flood in over a century in
 the summer of 1998 [del Ninno et al. (\citeyear{delninno2001}), del Ninno and
 Lundberg (\citeyear{delninno2005})].  The survey was fielded in three waves from
a sample of 757 households: round 1 in November, 1998; round
 2 in March--April, 1999; round 3 in November, 1999.  The
survey recorded all episodes of diarrhea for each child in the
household in the past six months or since the last interview by asking
the families at the time of each interview.  In addition,
 the survey recorded at each of the three interview times several
 important time-dependent covariates for the effect of diarrhea on
a child's future height: the child's current height and weight; the
amount of flooding in the child's home and village;  the
 household's economic and sanitation status.  In particular, the child's current height and weight
  are time-dependent confounders that satisfy conditions
 (c1) and (c2), making standard longitudinal data analysis methods
 biased [see Martorell and Ho (\citeyear{martorell1984}) and Moore et al. (\citeyear{moore2001}) for
 discussion of evidence for and reasons why current height and weight
 satisfy conditions (c1) and (c2)].  The time-dependent
confounders of current height and weight are available only at the time
of the interview,
 and changes in their value that might affect the exposure
of the child to the ``treatment'' of diarrhea, which varies in
continuous time, are not recorded in continuous time.
\end{example}

\begin{example}[{[The effect of AZT (Zidovudine) on CD4 counts]}]\label{exam2}
The Multicenter AIDS Cohort Study [MACS, Kaslow et al.
(\citeyear{kaslow1987})] has been used to study the effect of AZT on
CD4 counts [Hern\'{a}n, Brumback and Robins (\citeyear{hernan2002}),
Brumback et al. (\citeyear{brumack2004})]. Participants in the study
are asked to come semi-annually for visits at which they are asked to
complete a detailed interview, including a complete history of their
AZT use,
 as well as to take a physical examination.  Decisions on AZT
use are made by subjects and their physicians, and switches of
treatment might happen at any time between two visits.  These decisions
are based on the values of diagnostic variables, possibly including CD4
and CD8 counts, and the presence of certain symptoms. However, these
covariates are only measured by MACS at the time of visits; the values
of these covariates at the exact times that treatment decisions are
made between visits are not available.
\end{example}

\subsection{A model data generating process}\label{s:intro:dgp}
 In both the examples of AZT and diarrhea, the
exposure or treatment process happens continuously in time and a
complete record of the process is available, but the time-dependent
confounders are only observed at discrete times. There could be various
interpretations of the relationship between the data at the treatment
decision level and the data at the observational time level. To clarify
the problem of interest in this paper, we  consider a model data
generating process that satisfies all of the following assumptions:
\begin{longlist}
\item[(a1)] a patient takes a certain medicine under the advice of a
doctor;
\item[(a2)] a doctor continuously monitors and records a list
of health indicators
 of her patient and
decides the initiation and cessation of the medicine solely based on
current and historical records of these conditions, the historical use
of the medicine and possibly random factors unrelated to the patient's
health;
\item[(a3)] a third party organization asks a collection of
patients from various doctors to visit the organization's office
semi-annually; the organization measures the same list of health
indicators for the patients during their visits and asks the patients
to report the detailed history of the use of the medicine between two
visits;
\item[(a4)] we are only provided with the third party's data.
\end{longlist}
Note that in (a2), we assume the sequential randomization assumption
(A1) at the treatment decision level.

The AZT example can be approximated by the above data generating
process. In the AZT example, (a1) and (a2) approximately describe the
joint decision-making process by the patient and the doctor in the real
world.
 (a3) can be justified by reasonably assuming
that the staff at the MACS receive similar medical training and use
similar medical equipment as the patients' doctors. In the diarrhea
example, the patient's body, rather than a doctor, determines whether
the patient gets diarrhea. Assumption (a3), then, is saying that the
third party organization (the survey organization) collects enough
health data and that if all the histories of such health data are
available, the occurrence of diarrhea is conditionally independent of
the potential height.

\subsection{Difficulties posed by treatments varying in continuous
  time when covariates are observed only at discrete times}\label{s:intro:difficulty}

 Suppose our data are generated as in the previous section and we apply
discrete-time g-estimation at the discrete times at which the
time-dependent covariates are observed; we will denote these
observation times by $0,\ldots ,K$. In discrete-time g-estimation, we
are testing whether the observed treatment at time $t$ ($t=0,\ldots,K$)
is, conditional on the observed treatments at times $0,1,\ldots,t-1$
and observed covariates at times $0,\ldots,t$, independent of the
\textit{putative} potential outcomes at times $t+1,\ldots,K,$
calculated under the hypothesized treatment effect, where the putative
potential outcomes considered are what the subject's outcome would be
at times $t+1,\ldots,K$ if the subject never received treatment at any
time point. The difficulty with this procedure is that even if
sequential randomization holds when the measured confounders are
measured in continuous time [as is assumed in (a2)], it may not hold
when the measured confounders are measured only at discrete times. For
the discrete-time data, there can be \textit{unmeasured confounders}.
In the MACS example, the diagnostic measures at the time of AZT
initiation are missing unless the start of AZT initiation occurred
exactly at one of the discrete times that the covariates are observed;
the diagnostic measures at the initiation time are clearly important
confounders for the treatment status at the subsequent observational
time. In the diarrhea example, the nutrition status of the child before
the start of a diarrhea episode is missing unless the start of the
diarrhea episode occurred exactly at one of the discrete times that
covariates are observed; this nutrition status is also an important
confounder for the diarrhea status at the subsequent observational
time.
 Continuous-time sequential randomization does not, in
general, justify sequential randomization holding for the discrete-time
data, meaning that discrete-time g-estimation can produce inconsistent
estimates, even when continuous-time sequential randomization holds.

In this paper, we approach this problem from two perspectives. First,
we give conditions on the underlying continuous-time processes under
which discrete-time sequential randomization is implied, warranting the
use of discrete-time g-estimation.  Second, we propose a new estimation
method, called the \textit{controlling-the-future method}, that can
produce consistent estimates whenever discrete-time g-estimation is
consistent and can produce consistent estimates in some cases where
discrete-time g-estimation is inconsistent.

Our discussion focuses on a binary treatment and repeated continuous
outcomes. We also assume that the cumulative amount of treatment
between two visits is observed.  This is true for Examples \ref{exam1}
and \ref{exam2}, the AZT and diarrhea studies, respectively.  If
cumulative
 treatment is not observed,
there will often be a measurement error problem in the amount of
treatment, which is beyond the scope of this paper and an issue which
we are currently researching.

The organization of the paper is as follows: Section \ref{s:review}
reviews the standard discrete-time structural nested model and
g-estimation, describes a modified application when the underlying
process is in continuous time and proposes conditions on the
continuous-time processes when it works; Section \ref{s:new} describes
our controlling-the-future method; Section \ref{s:sim} presents a
simulation study; Section \ref{s:appl} provides an application to the
diarrhea study discussed in Example \ref{exam1};  Section \ref{s:concl}
concludes the paper.

\section{A modified g-estimation for discretely observed continuous-time processes}\label{s:review}
 In this section, we first review the discrete-time
structural nested model and the standard g-estimation, and
mathematically formalize the setting we described in Section
\ref{s:intro:dgp}. Then, with a slight modification and different
interpretation of notation, the g-estimation can be applied to the
discrete-time observations from the continuous-time model.
 We will show that under certain conditions, this estimation method is consistent.

\subsection{Review of discrete-time structural nested model and g-estimation}\label{s:review:dts}
 To reduce notation for the continuous-time
setting, we use a star superscript on every variable in this section.

Assuming that all variables can only change values at time $0, 1, 2,
\ldots, K$, we use $A^*_k$ to denote the binary treatment decision at
time $k$.
 Under the discrete-time setup, $A^*_k$
is assumed to be the constant level of treatment between time $k$ and
time $(k+1)$. We use $Y^{0*}_k$ to denote the baseline potential
outcome of the study at time $k$ if the subject does not receive any
treatment throughout the study and $Y^*_k$ to denote the actual outcome
at time $k$. In this paper, we assume that all $Y^{0*}_k$'s and
$Y^*_k$'s are continuous variables. Let $L^*_k$ be the vector of
covariates collected at time $k$. As a convention, $Y^*_k$ is included
in $L^*_k$.

We consider a simple deterministic model for the purposes of
illustration,
\begin{equation}\label{eq:ddc}
Y^*_k = Y^{0*}_k + \Psi \sum^{k-1}_{i=0} A^*_i,
\end{equation}
where $\Psi$ is the causal parameter of interest and can be interpreted
as the effect of one unit of the treatment on the outcome.

Model (\ref{eq:ddc}) is known as a \textit{rank-preserving} model
[Robins (\citeyear{robins1992})]. Under this model, for subjects $i$
and $j$ who have the same observed treatment history up to time $k$, if
we observe $Y_{k,i} < Y_{k,j}$, then we must have $Y^{0*}_{k,i} <
Y^{0*}_{k,j}$. It is \vspace*{1pt}also stronger than a more general rank-preserving
model since $Y^*_k$ depends deterministically only on $Y^{0*}_k$ and
the $A^*_i$'s.

Causal inference aims to estimate $\Psi$ from the observables, the
$A^*_k$'s and $L^*_k$'s. One way to achieve the identification of
$\Psi$ is to assume sequential randomization (A1). Given this notation
and model (\ref{eq:ddc}), a mathematical formulation of (A1) is
\begin{equation}\label{eq:seqrand_d}
    P(A^*_k|\bar{L}^*_k,\bar{A}^*_{k-1},\underline{Y}^{0*}_{k+})
=P(A^*_k|\bar{L}^*_k,\bar{A}^*_{k-1}),
\end{equation}
where $\bar{L}^*_k = (L_0, L_1, \ldots, L_k)$, $\bar{A}^*_{k-1} = (A_0,
A_1, \ldots, A_{k-1})$ and $\underline{Y}^{0*}_{k+}=(Y^{0*}_{k+1},
\break Y^{0*}_{k+2},\ldots, Y^{0*}_K)$.

For any hypothesized value of $\Psi$, we define a putative potential
outcome,
\[
Y^{0*}_k(\Psi) = Y^*_k - \Psi \sum^{k-1}_{i=0} A_i.
\]
Then, under (\ref{eq:ddc}) and (\ref{eq:seqrand_d}), the correct $\Psi$
should solve
\begin{equation}\label{eq:esteq_d}
  E[U(\Psi)] \equiv E\biggl\{\mathop{\sum_{k<m \leq K}}_{1 \leq i \leq N}[A^*_{i,k} -
      p_k(X^*_{i,k})]g(Y^{0*}_{i,m}(\Psi),X^*_{i,k})\biggr\} = 0,
\end{equation}
where $i$ is the index for each subject where there are $N$ subjects,
$X^*_{i,k} = (\bar{L}^*_{i,k}, \bar{A}^*_{i,k-1})$, $p_k(X^*_{i,k}) =
P(A^*_{i,k}=1 | X^*_{i,k}) $ is the propensity score for subject $i$ at
time $k$ and $g$ is any function. This estimating equation can be
generalized, with $g$ being a function of any number of future
$Y^{0*}_{i,m}(\Psi)$'s and $X^*_{i,k}$.

To estimate $\Psi$, we solve the empirical version of
(\ref{eq:esteq_d}):
\begin{equation} \label{eq:esteqdata_d}
    U(\Psi) \equiv \mathop{\sum_{k<m \leq K}}_{1 \leq i \leq N}[A^*_{i,k} -
    p_k(X^*_{i,k})]g(Y^{0*}_{i,m}(\Psi),X^*_{i,k})  = 0.
\end{equation}
If the true propensity score model is unknown and is parameterized as
$p_k(X^*_k, \beta)$, additional estimating equations are needed to
identify  $\beta$. For example, the following estimating equations
could be used:
\begin{equation}  \label{eq:esteqs_d}
    U(\Psi ,\beta ) = \mathop{\sum_{k<m \leq K}}_{1 \leq i \leq N}[A^*_{i,k} -
    p_k(X^*_{i,k})][g(Y^{0*}_{i,m}(\Psi),X^*_{i,k}),X^*_{i,k}]^T  = 0.
\end{equation}

The method is known as \textit{g-estimation}. The efficiency of the
estimate depends on the functional form of $g$. The optimal $g$
function that produces the most efficient estimation can be derived
[Robins (\citeyear{robins1992})]. The formulas for estimating the
covariance matrix of $(\hat{\Psi},\hat{\beta})$ are given in Appendix
\ref{appendixa}. A short discussion of the existence of the solution to
the estimating equation and identification can be found in Appendix
\ref{appendixb}.

\subsection{A continuous-time deterministic model and
  continuous-time sequential randomization}\label{s:review:model}
 We now extend the model in
Section \ref{s:review:dts} to a continuous-time model and define a
continuous-time version of the sequential randomization assumption (A1)
as a counterpart of (\ref{eq:seqrand_d}).

We now assume that the variables can change their values at any real
time between 0 and $K$. The model in Section \ref{s:review:dts} is then
extended as follows:

\begin{longlist}
\item[$\bullet$] $\{Y_t; 0 \leq t \leq K \}$ is the continuous-time,
  continuously-valued outcome process;
\item[$\bullet$] $\{L_t; 0 \leq t \leq K\}$ is the continuous-time
covariate
  process---it can be multidimensional and $Y_t$ is an element of
  $L_t$;
\item[$\bullet$] $\{A_t; 0 \leq t \leq K\}$ is the continuous-time
binary
  treatment process;
\item[$\bullet$] $\{Y^0_t; 0 \leq t \leq K \}$ is the continuous-time,
  continuously-valued potential outcome process if the subject does not
  receive any treatment from time 0 to time $K$---it can be thought
  of as the \textit{natural process} of the subject, free of
  treatment/intervention.
\end{longlist}

As a regularity condition, we further assume that all of the
continuous-time stochastic processes are c\`{a}dl\`{a}g processes
(i.e., continuous from the right, having limits from the left)
throughout this paper.

A natural extension of model (\ref{eq:ddc}) is
\begin{equation}\label{eq:dcc}
Y_t = Y^0_t + \Psi \int^t_0 A_s\,ds,
\end{equation}
where $\Psi$ is the causal parameter of interest. $\Psi$ can be
interpreted as the effect rate of the treatment on the outcome.

In this continuous-time model, a continuous-time version of the
sequential randomization assumption (A1) or, equivalently, assumption
(a2), can be formalized, although it does not have a simple form
similar to equation (\ref{eq:seqrand_d}). It was noted by Lok
(\citeyear{lok2008}) that a direct extension of the formula
(\ref{eq:seqrand_d}) involves ``conditioning null events on
  null events.''

Lok (\citeyear{lok2008}) formally defined continuous-time sequential
randomization when there is only one outcome at the end of the study.
We propose a similar definition for studies with repeated outcomes
under the deterministic model (\ref{eq:dcc}).

Let  $Z_t = (L_t, A_t, Y^0_t)$. Let $\sigma(Z_t)$ be the $\sigma$-field
generated \vspace*{1.5pt}by $Z_t$, that is, the smallest $\sigma$-field that makes
$Z_t$ measurable. Let $\sigma(\bar{Z}_t)$ be the $\sigma$-field
generated by $\bigcup_{u
  \leq t} \sigma(Z_u)$. Similarly,
$\sigma(\bar{Z}_t,\underline{Y}^0_{t+})$ is the $\sigma$-field
generated by $\sigma(\bar{Z}_t) \cup \sigma(\underline{Y}^0_{t+})$,
where $\sigma(\underline{Y}^0_{t+})$ is the $\sigma$-field generated by
$\bigcup_{u > t} \sigma(Y^0_u)$. By definition, the sequence of
$\sigma(\bar{Z}_t)$, $0 \leq t \leq K$, forms a filtration. The
sequence of $\sigma(\bar{Z}_t, \underline{Y}^0_{t+})$, $0 \leq t \leq
K$, also forms a filtration because $\sigma(\bar{Z}_t,
\underline{Y}^0_{t+}) \subset \sigma(\bar{Z}_s, \underline{Y}^0_{s+})$
for $t < s$ [note that this is true under the deterministic model
(\ref{eq:dcc}), but not in general].

Let $N_t$ be a counting process determined by $A_t$. It counts the
number of jumps in the $A_t$ process. Let $\lambda_t$ be a version of
the intensity process of $N_t$ with respect to $\sigma(\bar{Z}_t)$.
$M_t = N_t - \int^t_0 \lambda_s\,ds$  will be a martingale with respect
to $\sigma(\bar{Z}_t)$.

\begin{definition}\label{def:cti}
 With $N_t$ and $M_t$ defined as above, the
c\`{a}dl\`{a}g process $Z_t \equiv (L_t, A_t, Y^0_t)$, $0 \leq t \leq
K$, is said to satisfy the \textit{continuous-time sequential
randomization} assumption, or \textit{CTSR}, if $M_t$ is also a
martingale with respect to $\sigma(\bar{Z}_t,\underline{Y}^0_{t+})$.
Or, equivalently, there exists a $\lambda_t$ that is the intensity of
$N_t$, with respect to both the filtration of
$\sigma(\bar{Z}_t,\underline{Y}^0_{t+})$ and the filtration of
$\sigma(\bar{Z}_t)$.
\end{definition}

In this definition, given $A_0$, the counting process $\{N_t\}^T_0$
offers an alternative description of the treatment process
$\{A_t\}^T_0$. The intensity process $\lambda_t$, which models the
jumping rate of $N_t$, plays the same role as the propensity scores in
the discrete-time model, which models the switching of the treatment
process. Definition \ref{def:cti} formalizes  assumption (A1) in the
continuous-time model, by stating that $\lambda_t$ does not depend on
future potential outcomes.

The definition can be generalized if $A_t$ has more than two levels,
where $N_t$ can be a multivariate counting process, each element counts
a type of jump of the $A_t$ process and $\lambda_t$ is the multivariate
intensity process for $N_t$ under both the filtration of
$\sigma(\bar{Z}_t)$ and the filtration of
$\sigma(\bar{Z}_t,\underline{Y}^0_{t+})$; see Lok (\citeyear{lok2008}).

\subsection{A modified g-estimation}\label{s:review:g}
 In this paper, we assume that the continuous process
defined in Section \ref{s:review:model} can only be observed at integer
times, namely, times $0, 1, 2, \ldots, K$. We use the same  starred
notation as in Section \ref{s:review:dts}, but interpret instances of
this as discrete-time observations from the model in Section~\ref{s:review:model}.
Specifically:
\begin{longlist}
\item[$\bullet$] $\{A^*_k, k=0,1,2,\ldots,K\}$ denotes the set of
treatment
  assignments observable at times $0, 1, 2, \ldots, K$. We use
  $\bar{A}^*_k$ to denote the observed history of observed discrete-time
  treatment up to time $k$, that is, $(A^*_0$, $A^*_1$, \ldots, $A^*_k)$. Additionally,
  we use $\operatorname{cum}A^*_k = \int^{k-}_0 A_s\,ds$ to denote the cumulative amount
  of treatment up to time~$k$. Note that in the continuous-time model,
  $\operatorname{cum}A^*_k \neq \sum^{k-1}_{k'=0} A^*_{k'}$, as it would in discrete-time models.
  We let $\overline{\operatorname{cum}A^*_k} = (\operatorname{cum}A^*_1, \operatorname{cum}A^*_2, \ldots, \operatorname{cum}A^*_k)$.
  We note that, in practice, people sometimes use $\tilde{A}^*_k = \operatorname{cum}A^*_{k+1} - \operatorname{cum}A^*_k$
  as the treatment at time $k$ when applying discrete-time g-estimation to discrete-time
  observational data. Under deterministic models, such use of g-estimation usually
  requires stronger conditions than
  the conditions discussed in this paper. Throughout this paper,
   we define the treatment at time $k$ as $A^*_k$.

\item[$\bullet$] We define $L^*_k$, the observed covariates at
  time $k,$ to be $L_{k-}$, the left limit of $L$ at time $k$,
  following the convention that in the discrete model, people usually
  assume that the covariates are measured just before the treatment
  decision at time~$k$. $Y^*_k$ and $Y^{0*}_k$ are also \vspace*{1pt}defined as $Y_{k-}$ and
  $Y^0_{k-}$, respectively, following the same
  convention. $\bar{L}^*_{k}$ denotes\vspace*{1.5pt} $(L^*_0,
  L^*_1,\ldots,L^*_k)$, and $\bar{Y}^*_k$ and $\bar{Y}^{0*}_k$ are
  defined accordingly. $\underline{Y}^{0*}_{k+} = (Y^{0*}_{k+1},
  Y^{0*}_{k+2}, \ldots, Y^{0*}_K)$.
\end{longlist}

With this notation and in the spirit of g-estimation, which controls
all observed history in the propensity score model for the treatment,
we  propose the following working estimating equation:
\begin{equation}
    U(\Psi) \equiv \mathop{\sum_{k<m \leq K}}_{1 \leq i \leq N}[A^*_{i,k} -
    p_k(X^*_{i,k})]g(Y^{0*}_{i,m}(\Psi),X^*_{i,k}) = 0,
  \label{eq:esteqdata}
\end{equation}
where $X^*_{i,k}$ is the collection of $\bar{L}^*_{i,k},
\bar{A}^*_{i,k-1}$ and $\overline{\operatorname{cum}A^*_{i,k}}$,
$p_k(X^*_{i,k}) = P(A^*_{i,k}=1 | X^*_{i,k})$ and\vspace*{1pt} $Y^{0*}_{i,m}(\Psi) =
Y^*_{i,m} - \Psi \operatorname{cum}A^*_{i,k}$.

In practice, $p_k(X^*_{i,k})$ is unknown and has to be\vspace*{1.5pt} parameterized as
$p_k(X^*_{i,k};$ $\beta)$, and we use different functions $g$ to
identify all of the parameters, as in Section~\ref{s:review:dts}. The
covariance matrix of estimated parameters can be estimated as in
Appendix~\ref{appendixa}. A discussion of the existence of a solution
and identification can be found in Appendix \ref{appendixb}.

The estimating equation has the same form as (\ref{eq:esteqdata_d}),
except for two important differences. First, the propensity score model
in this section conditions on the additional
$\overline{\operatorname{cum}A^*_{i,k}}$. In the discrete-time model of
Section \ref{s:review:dts}, $\overline{\operatorname{cum}A^*_{i,k}}$
would be a transformed version of $\bar{A}^*_{i,k-1}$ and was redundant
information. However, with continuous-time underlying processes,
$\overline{\operatorname{cum}A^*_{i,k}}$ provides new information on
the treatment history. Second, the putative potential outcome
$Y^{0*}_{i,m}(\Psi)$ is calculated by subtracting the $\operatorname{cum}A^*_{i,k}$
from $Y^*_{i,m}$, instead of $\sum^{k-1}_{l=0} A^*_{i,l}$. We will
later refer to the g-estimation in this section as the \textit{modified
g-estimation} (although it is in the true spirit of g-estimation). The
justification and limitation of using the modified g-estimation will be
discussed in Section \ref{s:just}.

We refer to the g-estimation in Section \ref{s:review:dts} as
\textit{naive g-estimation} when it is applied to data from a
continuous-time model. When the data come from a continuous-time model,
the naive g-estimation can be severely biased, as we will show in our
simulation study and the diarrhea application. One source of bias is a
measurement error problem, $\sum^{k-1}_{l=0}A_{i,l}^*$ is not the
correct measure of the treatment; another source of bias is that the
important information $\overline{\operatorname{cum}A_{i,k}}$ is not
conditioned on in the propensity score. Although we would not expect
researchers to use naive g-estimation when the true cumulative
treatments are available, we  present the simulation and real
application results using this method as a reference to show how
severely biased the estimates would be had we not known the true
cumulative treatments and  the measurement error problem had dominated.

\subsection{Justification of the modified g-estimation}\label{s:just}
 Given discrete-time observational data from
continuous-time underlying processes, solving equation
(\ref{eq:esteqdata}) provides an estimate for $\Psi$. For this $\Psi$
estimate to be consistent, an analog to  condition (\ref{eq:seqrand_d})
is needed:
\begin{equation} \label{eq:seqrand}
    P(A^*_k|\bar{L}^*_k,\bar{A}^*_{k-1},\overline{\operatorname{cum}A^*_k}, \underline{Y}^{*0}_{k+})
=P(A^*_k|\bar{L}^*_k,\bar{A}^*_{k-1},\overline{\operatorname{cum}A^*_k}).
\end{equation}

Condition (\ref{eq:seqrand}) is a requirement on variables at
observational time points. Its validity for a given study relies on how
the data are collected, in addition to the underlying continuous-time
data generating process. It is not clear, without conditions on the
underlying continuous-time data generating process, how one would go
about collecting data in a way such that (\ref{eq:seqrand}) would hold
while the standard ignorability (\ref{eq:seqrand_d}) is not true. Here,
we will seek conditions at the continuous-time process level that imply
condition (\ref{eq:seqrand}) and hence justify the estimating equation~(\ref{eq:esteqdata}).
In particular, we consider two such conditions.

\subsubsection{Sequential randomization at any finite subset of time
  points}\label{s:just:subset}
Recall the data generating process described in
Section \ref{s:intro:dgp}.  The third party organization periodically
(e.g., semi-annually) collects the health data and treatment records of
the patients.  Suppose that a researcher thinks (\ref{eq:seqrand})
holds for the time points at which the third party organization
collects these data.  If the time points have not been chosen in a
special way to make (\ref{eq:seqrand}) hold, then the researcher will
often be willing to make the stronger assumption that
(\ref{eq:seqrand}) would hold for any finite subset of time points at
which the third party organization chose to collect data.  For example,
for the diarrhea study, the survey was actually conducted in November,
1998, March--April, 1999 and November, 1999.  If a researcher thought
(\ref{eq:seqrand}) held for these three time points, then she might be willing to
assume that (\ref{eq:seqrand}) should also hold if  the survey was
instead conducted in December, 1998, February, 1999, May, 1999 and
October, 1999.

Before formalizing the researcher's assumption on any finite subset of
time points, we make the following
observation.
\begin{proposition}\label{prop:y0}
 Under the deterministic model assumption
(\ref{eq:dcc}), the propensity score has the following property:
\begin{equation}\label{eq:seqrand_y0}
P(A^*_k=1|\bar{L}^*_k,\bar{A}^*_{k-1},\overline{\operatorname{cum}A^*_{k}})
= P(A^*_k=1 | \bar{L}^*_k, \bar{A}^*_{k-1}, \bar{Y}^{0*}_{k} ).
\end{equation}
\end{proposition}

\begin{pf} Under the deterministic assumption
(\ref{eq:dcc}) and the correct $\Psi$, $(\bar{L}^*_k,$
$\bar{A}^*_{k-1},$ $\overline{\operatorname{cum}A^*_k})$ is  a
one-to-one transformation of $(\bar{L}^*_k, \bar{A}^*_{k-1},$
$\bar{Y}^{0*}_{k})$.
\end{pf}

Using Proposition \ref{prop:y0}, we state the sequential randomization
assumption at any finite subset of time points as follows.
\begin{definition}\label{def:fti}
 A c\`{a}dl\`{a}g process $Z_t \equiv (L_t, A_t,
Y^0_t)$, $0 \leq t \leq K$, is said to satisfy the \textit{finite-time
sequential randomization} assumption, or \textit{FTSR}, if, for any
finite subset of time points, $0 \leq t_1 < t_2 < \cdots < t_n <
t_{n+1} < \cdots < t_{n+l} \leq K$, we have
\begin{equation}\label{eq:fti}
  P(A_{t_n} | \bar{L}_{t_n-}, \bar{A}_{t_{n-1}}, \bar{Y}^0_{t_n-},
  \underline{Y}^0_{t_n +}) = P(A_{t_n} | \bar{L}_{t_n-},
  \bar{A}_{t_{n-1}}, \bar{Y}^0_{t_n-}),
\end{equation}
where $\bar{L}_{t_n-} = (L_{t_1-},L_{t_2-},\ldots,L_{t_n - })$,
$\bar{A}_{t_{n-1}} = (A_{t_1},A_{t_2},\ldots,A_{t_{n-1}})$,
$\bar{Y}^0_{t_n-} = (Y^0_{t_1-}, Y^0_{t_2-},\ldots,Y^0_{t_n - })$
and $\underline{Y}^0_{t_n+} =
(Y^0_{t_{n+1}-},Y^0_{t_{n+2}-},\ldots,Y^0_{t_{n+l}-})$.
\end{definition}

It should be noted that for the conditional densities in
(\ref{eq:seqrand_y0}) and (\ref{eq:fti}), and the conditional densities
in the following sections, we always choose the version that is the
ratio of joint density to marginal density.

The finite-time sequential randomization assumption clearly implies
condition (\ref{eq:seqrand_d}) and thus justifies the modified
g-estimation equation (\ref{eq:esteqdata}). We have also proven a
result that shows the relationship between the FTSR assumption and the
CTSR assumption.

\begin{theorem}\label{thm:main}
 If a continuous-time c\`{a}dl\`{a}g process $Z_t$
satisfies finite-time sequential randomization, then, under some
regularity conditions, it will also satisfy continuous-time sequential
randomization.
\end{theorem}

\begin{pf}
See Appendix \ref{appendixc}. The regularity conditions are also stated
in Appendix \ref{appendixc}.
\end{pf}

The result of Theorem \ref{thm:main} is natural. As mentioned in
Section \ref{s:intro:difficulty}, the continuous-time sequential
randomization does not imply FTSR because, in discrete-time
observations, we do not have the full continuous-time history to
control. To compensate for the incomplete data problem, some stronger
assumption on the continuous-time processes must be made if
identification is to be achieved.

\subsubsection{A Markovian condition}\label{s:just:cond}
 Given the finite-time sequential randomization
assumption described above, two important questions arise. First,
Theorem~\ref{thm:main} shows that the FTSR assumption is stronger than
the continuous-time sequential randomization assumption. It is natural
to ask how much stronger it is than the CTSR assumption. Second, the
FTSR assumption, unlike the CTSR assumption (A1), is not an assumption
on the data generating process itself and so it is not clear how to
incorporate domain knowledge about the data generating process to
justify it. Is there a condition at the data generating process level
which will be more helpful in deciding whether g-estimation is valid?

We partially answer both questions in the following theorem.

\begin{theorem}\label{thm:cond}
 Assuming that the process $(Y^0_t, L_t, A_t)$
satisfies the \textit{continuous-time sequential randomization}
assumption, and that the process $(Y^0_{t-}, L_{t-}, A_t)$ is Markovian,
for any time $t$ and $t+s$, $s>0$, we have
\begin{equation}
P(A_t | L_{t-}, Y^0_{t-}, \underline{Y}^0_{t +}) = P(A_t | L_{t-} ,
Y^0_{t-}), \label{eq:markovcon}
\end{equation}
which implies the \textit{finite-time sequential randomization}
assumption. Here, $\underline{Y}^0_{t+} =
(Y^0_{t_{1}-},Y^0_{t_{2}-},\ldots,Y^0_{t_{n}-})$ and $t < t_1 < t_2 <
\cdots < t_n$.
\end{theorem}

\begin{pf} The proof can be found in the Appendix \ref{appendixd}.
\end{pf}

The theory states that the Markov condition and the CTSR assumption
together imply the FTSR condition. Therefore, they imply condition
(\ref{eq:seqrand_d}) and thus justify the modified g-estimation
equation (\ref{eq:esteqdata}).

We make the following comments on the theorem.
\begin{longlist}
\item[$\bullet$] The theorem partially answers our first question---the
FTSR assumption is stronger than
    the CTSR assumption, but the gap between the two assumptions is less than a Markovian assumption.
    The result is not surprising since, with missing covariates between observational time points, we
    would hope that the variables at the observational time points well summarize the missing information.
     The Markovian assumption guarantees that variables at an observational time point summarize all information prior to that time point.
\item[$\bullet$] The theorem also partially answers our second
question. The CTSR assumption is usually
    justified by domain knowledge of how treatments are decided. Theorem \ref{thm:cond}
     suggests that the researchers could further look for biological evidence that the process
      is Markovian to validate the use of g-estimation. The Markovian assumption can also be tested.
    One could first use the modified g-estimation to estimate the causal parameter,
    construct the $Y^0$ process at the observational time points and then test whether the full
     observational data
    of $A, L, Y^0$ come from a Markov process. A strict test of whether the discretely observed
     longitudinal data come from a continuous-time (usually nonstationary) Markov process could
      be difficult and is beyond the scope of this paper. As a starting point, we suggest Singer's trace inequalities [Singer (\citeyear{singer1981})] as a
    criterion to test for the Markovian property. A weaker test for the Markovian property is to
     test conditional independence of past observed values and future observed values conditioning on current observed values.
\item[$\bullet$] In the theorem, equation (\ref{eq:markovcon}) looks
like an even stronger version of the continuous-time sequential
randomization assumption---the treatment decision seems to be based
only on current covariates and current potential outcomes. One could,
of course, directly assume this stronger version of randomization and
apply g-estimation. However, Theorem \ref{thm:cond} is more useful
since we are assuming a weaker untestable CTSR assumption and a
Markovian assumption that is testable in principle.
\item[$\bullet$]
The theorem suggests that it is sufficient to control for current
covariates and current potential outcomes for g-estimation to be
consistent. In practice, we advise controlling for necessary past
covariates and treatment history. The estimate would still be
consistent if the Markovian assumption were true and it might reduce
bias when the Markovian assumption was not true. As a result, we do
control for previous covariates and treatments in our simulation and
application to the diarrhea data.
\item[$\bullet$] It is worth noting that the
labeling of time is arbitrary. In practice, researchers can label
whatever they have controlled for in their propensity score as the
``current'' covariates, which could include covariates and treatments
that are measured or assigned previously. In this case, the dimension
of the process that needs to be tested for the Markovian property
should also be expanded to include older covariates and treatments.
\item[$\bullet$] Finally, we note that a discrete-time version of the theorem is
implied by Corollary 4.2 of Robins (\citeyear{robins1997}) if we set, in his notation,
$U_{ak}$ to be the covariates between two observational time points and
$U_{bk}$ to be the null set.
\end{longlist}

\begin{figure}
  \begin{tabular}{c}

\includegraphics{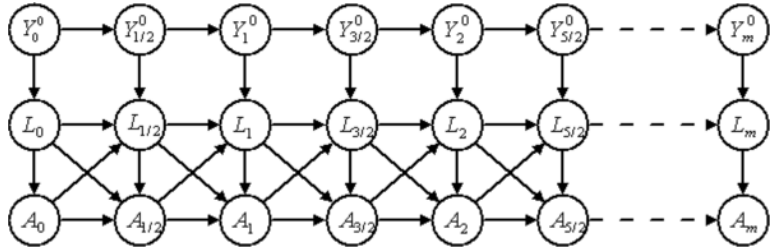}
\\ [-1pt]
 (a) DAG of a Markovian process.\\ [9pt]

\includegraphics{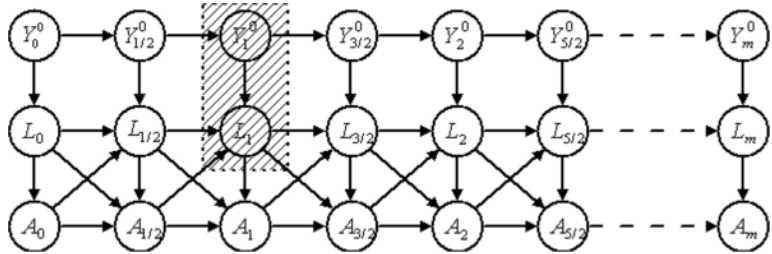}
\\ [-1pt]
 (b) Verification of equation (\protect\ref{eq:seqrand_y0}).
\end{tabular}
\caption{Directed acyclic graph.} \label{fig:DAG}
\end{figure}

As a discretized example, we illustrate the idea of Theorem
\ref{thm:cond} by a directed acyclic graph
  (DAG) in part (a) of Figure
  \ref{fig:DAG}, which assumes that all variables can only change
  values at time points  $0,1/2,1,3/2,2,\ldots,m$. Note
  that we do not distinguish the left limits of variables and the variables
  themselves in all DAGs of this paper, for reasons discussed in Appendix \ref{appendixc}.
  We also assume that the process can only be
  observed at times $0, 1, 2, \ldots,m$. It is easy to verify that the DAG satisfies
   sequential randomization at the $0,1/2,1,3/2,2,\ldots,m$ time level. The DAG is also Markovian in time. For example,
  if we control $A_1, L_1, Y^0_1$, any variable prior to time 1 will be d-separated from any other variable after time 1.

  Part (b) of Figure \ref{fig:DAG} verifies that $A_1$ is d-separated from
  $Y^0_m$, $m > 1$ by the shaded variables, namely, $L_1$ and
  $Y^0_1$, as is implied by equation (\ref{eq:seqrand_y0}). By
  Theorem \ref{thm:cond}, the modified g-estimation works for data
  observed at the integer times if they are generated by the model
  defined by this DAG.

It is true that the Markovian condition that justifies the g-estimation
equation (\ref{eq:esteqdata}) is restrictive, as will be discussed in
the following section. However, our simulation study shows that
g-estimation has some level of robustness when the Markovian assumption
is not seriously violated.

\section{The controlling-the-future method}\label{s:new}

In this section, we consider situations in which the
 observational time sequential randomization fails and seek methods that are
more robust to this failure than the modified g-estimation given in
Section \ref{s:review:g}. The method we are going to introduce was
proposed in Joffe and Robins (\citeyear{joffe2009}), which deals with a more general
case of
 the existence of unmeasured confounders. It can be applied to deal with unmeasured
confounders coming from either a subset of contemporaneous covariates
or a subset of covariates that represent past time, the latter case
being of interest for this paper. The method, which we will refer to as
the controlling-the-future
 method (the reason for the name will become clear later on), gives consistent estimates when
g-estimation is consistent and  produces consistent estimates in some
cases even when g-estimation is severely inconsistent.

In what follows, we will first describe an illustrative application of
the controlling-the-future method and then discuss its relationship
with our framework of g-estimation in continuous-time processes with
covariates observed at discrete times.

\subsection{Modified assumption and estimation of parameters}\label{s:new:mod}

We assume the same continuous-time model as in Section
\ref{s:review:model}. Following Joffe and Robins (\citeyear{joffe2009}), we consider a
revised sequential randomization assumption on variables at the
observational time points
\begin{equation}\label{eq:joffe}
P(A^*_k|\bar{L}^*_k,\bar{A}^*_{k-1},\overline{\operatorname{cum}A^*_{k}},
\underline{Y}^{0*}_{k+})=P(A^*_k|\bar{L}^*_k,
\bar{A}^*_{k-1},\overline{\operatorname{cum}A^*_{k}},Y^{0*}_{k+1}).
\end{equation}

This assumption relaxes (\ref{eq:seqrand}). At each time point,
conditioning on previous observed history, the treatment can depend on
future potential outcomes, but only on the next period's potential
outcome. In Joffe and Robins' extended formulation, this can be further
relaxed to allow for dependence on more than one period of future
potential outcomes, as well as other forms of dependence on the
potential outcomes.

If the revised assumption (\ref{eq:joffe}) is true, then we obtain a
similar
 estimating equation as (\ref{eq:esteqdata}).
 For each putative $\Psi$, we map $Y^*_k$ to
    \[
    Y^{0*}_k(\Psi) = Y^*_k - \Psi \operatorname{cum}A^*_{k},
    \]
the potential outcome if the subject never received any treatment under
the hypothesized treatment effect $\Psi$.

Define the putative propensity score as
\begin{equation}\label{eq:ppscore}
p_k(\Psi) \equiv
P(A^*_k=1|\bar{L}^*_k,\bar{A}^*_{k-1},\overline{\operatorname{cum}A^*_{k}},Y^{0*}_{k+1}(\Psi)).
\end{equation}

Under assumption (\ref{eq:joffe}), the correct $\Psi$ should solve
\begin{equation}\label{eq:joffeesteq}
\qquad U(\Psi) = E\biggl\{\mathop{\sum_{1 \leq i \leq n}}_{k+1<m \leq K
}[A^*_{i,k}-p_{i,k}(\Psi)]g(Y^{0*}_{i,m}(\Psi),
X^*_{i,k},h_{i,k}(\Psi))\biggr\} = 0,
\end{equation}
where $X^*_{i,k} =
(\bar{L}^*_{i,k},\bar{A}^*_{i,k-1},\overline{\operatorname{cum}A^*_{i,k}})$,
$h_{i,k}(\Psi) = Y^{0*}_{i,k+1}(\Psi)$ and $g$\vspace*{1.5pt} is any function and can
be generalized to functions of $X^*_{i,k}$, $h_{i,k}(\Psi)$ and any
number\vspace*{-2pt} of future potential outcomes that are later than time $k+1$, for
example, $g(Y^{0*}_{i,k+2}$$(\Psi),$ $Y^{0*}_{i,k+3}(\Psi),$
$X^*_{i,k},h_{i,k}(\Psi))$. In most real applications, the model\vspace*{1.5pt} for
$p_k(\Psi)=E[A^*_k|X^*_k,h_k(\Psi)]$ is unknown and is usually
estimated by a parametric model,
\[
p_{i,k}(\Psi;\beta_X,\beta_h) =
E[A_{i,k}|X^*_{i,k},h_{i,k}(\Psi);\beta_X,\beta_h].
\]

We can solve the following set of estimating equations to obtain the
estimates of $\Psi$, $\beta_X$ and $\beta_h$:
\begin{eqnarray}\label{eq:joffeesteqs}
U(\Psi,\beta_X,\beta_h)&=& \mathop{\sum_{1 \leq i \leq n}}_{ k+1<m \leq K }
\bigl(A^*_{i,k}-p_{i,k}(\Psi;\beta_X,\beta_h)\bigr)\nonumber\\
&&\hphantom{\mathop{\sum_{1 \leq i \leq n}}_{ k+1<m \leq K }}{}\times[g(Y^{0*}_{i,m}(\Psi),X^*_{i,k},h_{i,k}(\Psi)),X^*_{i,k},h_{i,k}(\Psi)]^T \\
&=& 0. \nonumber
\end{eqnarray}
The estimation of the covariance matrix of $\Psi$, $\beta_X$ and
$\beta_h$ is similar to the usual standard g-estimation, which is
described in Appendix \ref{appendixa}.

Two important features of estimating equation (\ref{eq:joffeesteqs})
distinguish it from estimating equation (\ref{eq:esteqdata}). First, in
(\ref{eq:joffeesteqs}), there is a common parameter $\Psi$ in both
$p_k$'s model and $Y^{0*}_{m}(\Psi)$, caused by the fact that the
treatment depends on a future potential outcome. Second, in
(\ref{eq:joffeesteqs}), the sum over $m$ and $k$ is restricted to
$m>k+1$,
 while in (\ref{eq:esteqdata}), we only need $m>k$. If
we use $m=k+1$ in (\ref{eq:joffeesteqs}),
$E\{[A^*_{i,k}-p_{i,k}(\Psi)]g(Y^{0*}_{i,k+1}(\Psi),
X^*_{i,k},h_{i,k}(\Psi))\}=0$ usually \vspace*{1pt}does not lead to the
identification of $\Psi$, unless certain functional forms of the
propensity score model are assumed to be true [see Joffe and Robins
(\citeyear{joffe2009})].

\subsection{The controlling-the-future method and the Markovian condition}\label{s:new:intui}

Joffe and Robins' revised assumption (\ref{eq:joffe}) is  an assumption
on the discrete-time observational data. It relaxes the observational
 time sequential randomization (\ref{eq:seqrand}) because
 (\ref{eq:seqrand}) always implies (\ref{eq:joffe}). At the
continuous-time data generating level, (\ref{eq:joffe}) allows less
stringent underlying stochastic processes  than the Markovian process
in Theorem \ref{thm:cond}.

In particular, we identify two important scenarios where the relaxation
happens. One scenario is to allow for more direct temporal dependence
for the $Y^0$ process, which we will refer to as the
\textit{non-Markovian-$Y^0$} case. The other scenario is to allow
colliders in $L$, which we will refer to as the
\textit{leading-indicator-in-$L$} case. We illustrate both cases by
modifying the directed acyclic graph (DAG) example in
Figure \ref{fig:DAG}.

\begin{figure}
\begin{tabular}{c}

\includegraphics{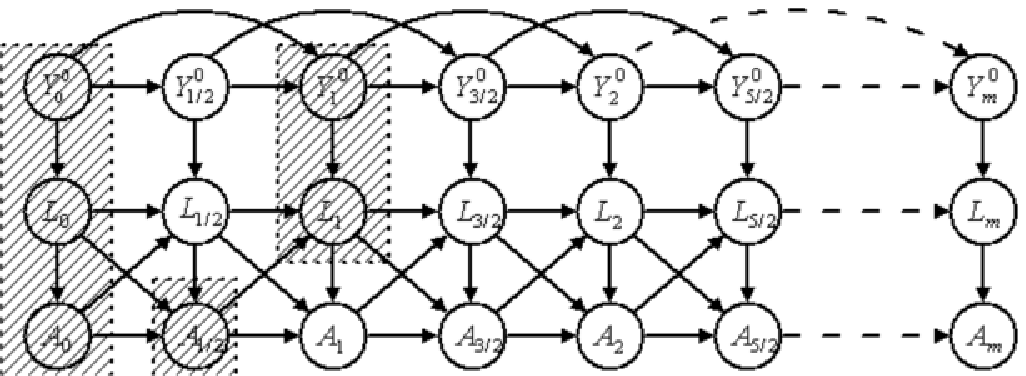}
\\ [-2pt]
  (a) No control for future $Y^0_t$.\\ [9pt]

\includegraphics{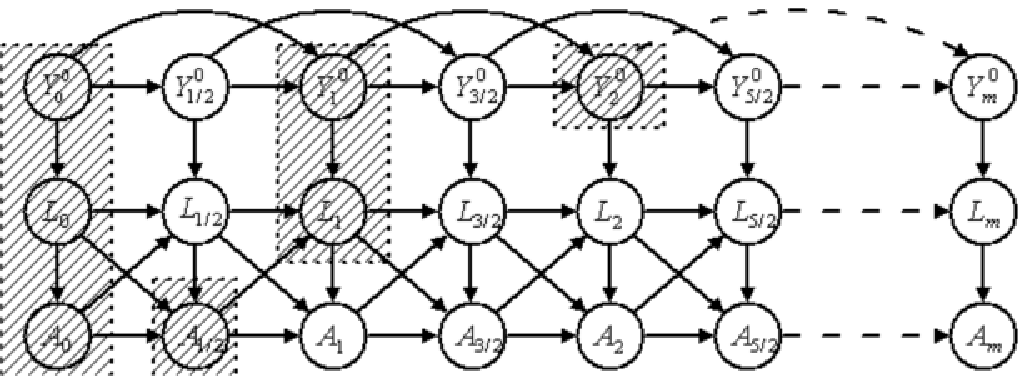}
\\ [-2pt]
  (b) Control for future $Y^0_t$.
\end{tabular}
\caption{Directed acyclic graph with non-Markovian $Y^0_t$.}\label{fig:nonmarkov}
\end{figure}

\begin{The non-Markovian-$Y^0$ case*}
Assume, for example, our data is generated from the DAG in
Figure \ref{fig:nonmarkov}, where we allow the dependence of $Y^0_2$ on
$Y^0_1$, even if $Y^0_{3/2}$ is controlled.  In part (a) of
Figure \ref{fig:nonmarkov}, we control\vspace*{1.5pt} for observed covariates ($L_0,
L_1$), treatment ($A_0, A_{1/2}$) and current and historical potential
outcome ($Y^0_0, Y^0_1$) for treatment at time 1 ($A_1$), that is, we
have controlled for all historically observed covariates, treatment and
cumulative treatment as suggested in the comments accompanying Theorem
\ref{thm:cond}. In this case, the modified g-estimation fails because
the paths\vspace*{1.5pt} like $A_1 \leftarrow L_{1/2} \leftarrow Y^0_{1/2} \rightarrow
Y^0_{3/2} \rightarrow Y^0_2 \rightarrow \cdots\rightarrow Y^0_m$ are not
blocked by the shaded variables. In part (b) of Figure
\ref{fig:nonmarkov}, we control for the additional $Y^0_2$. $A_1$ is
not completely blocked from $Y^0_m$, but some paths that are not
blocked in part (a) are now blocked, for example, the path of $A_1
\leftarrow
 L_{1/2} \leftarrow Y^0_{1/2} \rightarrow
 Y^0_{3/2} \rightarrow Y^0_2 \rightarrow Y^0_{5/2} \rightarrow
 \cdots \rightarrow
 Y^0_m$. Also, no additional paths are opened by conditioning on $Y_2 ^0$.
 We would usually expect that the correlation between $A_1$ and $Y^0_m$ is
 weakened.
Under the framework of Joffe and Robins (\citeyear{joffe2009}), we can control for more
than
 one period of future potential outcomes and expect to further weaken the
 correlation between $A_1$ and $Y^0_m$. A modification of assumption (\ref{eq:joffe}) that conditions on more future potential outcomes may be
 approximately true.

The scenario relates to real-world problems. For instance, in the
diarrhea example, $Y^0_t$ is the natural height growth of a child
without any occurrence of diarrhea. Height in the next month not only
depends on the current month's height, but also depends on the previous
month's height: the complete historical growth curve of the child
provides information on genetics and nutritional status, and provides
information about future natural height beyond that of
 current natural height alone. Therefore, the
potential height process for the child is not Markovian. [For a formal
argument why children's height growth is not Markovian, see Gasser et
al. (\citeyear{gasser1984}).] By the reasoning employed above, g-estimation fails.
However, if we assume that the delayed dependence of natural height
wanes after a period of time (as in Figure \ref{fig:nonmarkov}),
controlling for the next period potential height in the propensity
score model might weaken the relationship between current diarrhea
exposure and future potential height later than the next period and the
assumptions of the controlling-the-future method might hold
approximately.
\end{The non-Markovian-$Y^0$ case*}

\begin{The leading-indicator-in-$L$ case*}
In Figure \ref{fig:DAG}, we do not allow any arrows from future $Y^0$ to
previous $L$, which means that among all measures of the subject, there
are no elements in $L$ that contain any leading information about
future $Y^0$. This means that $Y^0$ is a measure that is ahead of all
other measures,  by which we mean that, for example, $L_2 \perp Y^0_m |
\bar{Y}^0_2, \bar{A}_{2-}, \bar{L}_{2-} m>2$. This is not realistic in
many real-world problems. In the example of the effect of the diarrhea
on height, weight is an important covariate. While both height and
weight reflect the nutritional status of a child,
 malnutrition usually affects weight more quickly than
height, that is, the weight contains leading information for the
natural height of the child.
 Figure \ref{fig:DAG} is thus not an appropriate model
for studying the effect of  diarrhea on height.
\begin{figure}
  \begin{tabular}{c}

\includegraphics{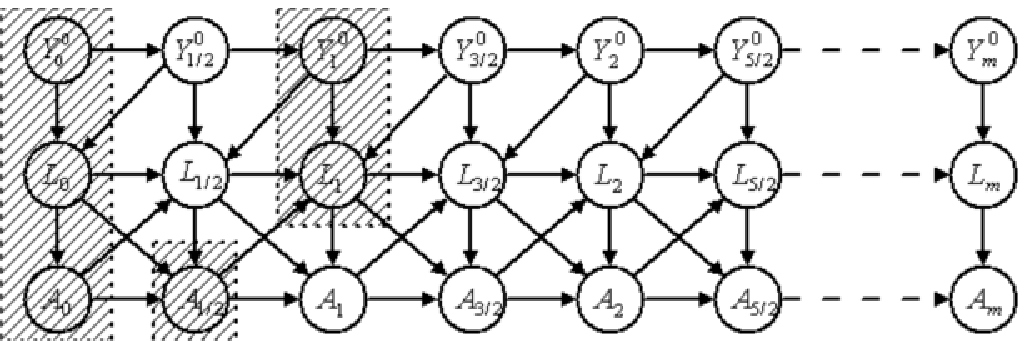}
\\ [-2pt]
   (a) No control for future $Y^0_t$.\\ [9pt]

\includegraphics{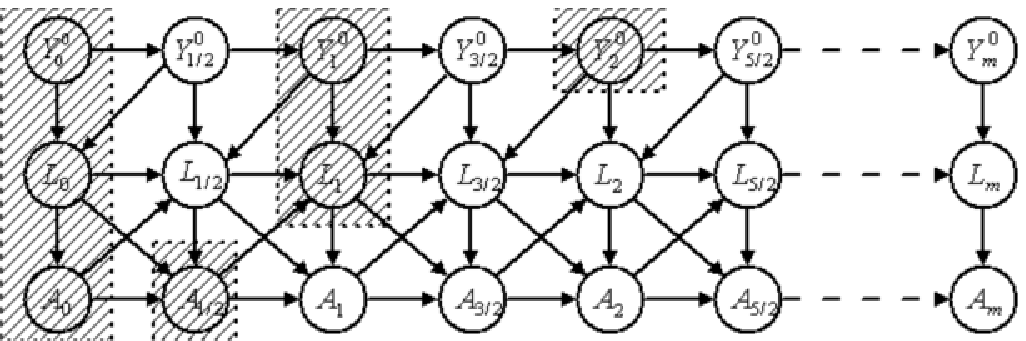}
\\ [-2pt]
   (b) Control for future $Y^0_t$.\\
\end{tabular}
\caption{Directed acyclic graph with leading indicator in $L_t$.}\label{fig:leadinginfo}
\end{figure}

In Figure \ref{fig:leadinginfo}, we allow arrows from $Y^0_{1/2}$ to
$L_0$, from $Y^0_1$ to $L_{1/2}$ and so on, which  assumes that $L$
contains leading indicators of $Y^0$, but the leading indicators are
only ahead of $Y^0$ for less than one unit of time. Part (a) of
Figure~\ref{fig:leadinginfo} shows that controlling for history of
covariates, treatment and potential outcomes does not block $A_1$ from
$Y^0_m$. On the path of $A_1 \leftarrow  L_{1/2} \rightarrow  L_1
\leftarrow  Y^0_{3/2} \rightarrow Y^0_2 \rightarrow  Y^0_{5/2}
\rightarrow  \cdots \rightarrow Y^0_m$, $L_1$ is a controlled collider.
However, in part (b), if we do control for $Y^0_2$ additionally, the
same path will be blocked. In general, if we assume that there exist
leading indicators in covariates and that the leading indicators are
not ahead of potential outcomes for more than one time unit,
g-estimation will fail, but the controlling-the-future method will
produce consistent estimates.

The fact that the controlling-the-future method can work in the leading
information scenario can also be related to the discussion of
 Section 3.6 of Rosenbaum (\citeyear{rosenbaum1984}). The main reason for g-estimation's
 failure in the DAG example is that $L_{1/2}$ is not observable and
 cannot be controlled. If $L_{1/2}$ is observed, it is easy to verify
that the DAG in Figure \ref{fig:leadinginfo} satisfies sequential
 randomization on the finest time grid.
 The idea behind
the controlling-the-future method is to condition on a ``surrogate'' for
$L_{1/2}$. The surrogate should satisfy the property that $Y_m^0$ is
independent of the unobserved $L_{1/2}$ given the surrogate and other
observed covariates [similar to formula 3.17 in Rosenbaum (\citeyear{rosenbaum1984})].  In
the leading information case, when  $m>k+1$ and we have covariates
$\bar{L}_k$ that are only ahead of the potential outcome until time at
most $k+1$, the future potential outcome $Y_{k+1}^0$ is a surrogate. It
is easy to check that in Figure \ref{fig:leadinginfo}, $L_{1/2}$ is
independent of $Y^0_m$, given $Y^0_2$, $L_1$, $A_0$ and $\operatorname{cum}A_1$
(equivalently, $Y^0_1$).

It is worth noting that we do not need to control for anything except
$Y_2^0$ in Figure~\ref{fig:leadinginfo} in order to get a consistent
estimate. It is possible to construct more complicated DAGs in which
controlling for additional past and current covariates is necessary,
which involves more model specifications for the relationships among
different covariates and deviates from the main point of this paper.
\end{The leading-indicator-in-$L$ case*}

In Section \ref{s:sim}, we will simulate data in cases of
non-Markovian-$Y^0_t$ and leading-indicator-in-$L_t$, respectively, and
show that the controlling-the-future method does produce better
estimates than  g-estimation.
 However, it is worth noting that when the modified g-estimation in
 Section \ref{s:review:g} is consistent, the controlling-the-future estimation is usually considerably less efficient.
 This is because condition (\ref{eq:joffe}) is less stringent than (\ref{eq:seqrand}).
  The semiparametric model under (\ref{eq:seqrand}) is a submodel of
   the semiparametric model defined by (\ref{eq:joffe}).
   The latter will have a larger semiparametric efficiency bound
    than the former. Theoretically, the most efficient g-estimation will
     be more efficient than the most efficient controlling-the-future estimation
      if the g-estimation is valid. In practice, even if we are not using the most
       efficient estimators, controlling-the-future estimation usually estimates more parameters,
        for example, coefficients for $h_{i,k}(\Psi)$ in the propensity model, and thus is less
         efficient. For a formal discussion, see Tsiatis (\citeyear{tsiatis2006}).

\section{Simulation study}\label{s:sim}

We set up a simple continuous-time model that satisfies sequential
ignorability in continuous time, and simulate and record discrete-time
data from variations of the simple model. We estimate causal parameters
from both the modified g-estimation and the controlling-the-future
estimation. We also present the estimates from naive g-estimation in
Section \ref{s:review:dts}, where we ignore the continuous-time
information of the treatment processes, as a way to show the severity
of the bias in the presence of the measurement error problem. The
results support the discussions in Sections \ref{s:just} and
\ref{s:new}.

In the simulation models below, M1 satisfies the Markovian condition in
Theorem \ref{thm:cond}. It also serves as a proof that there exist
processes satisfying the conditions of Theorem \ref{thm:cond}.

\subsection{The simulation models}\label{s:sim:mod}

We first consider a continuous-time Markov model which satisfies the
CTSR assumption.
\begin{longlist}
\item[$\bullet$] $Y^0_t$ is the potential outcome process if the patient is
  not receiving any treatment. We assume that
  \[
  Y^0_t = g(V,t) + e_t,
  \]
  where $g(V,t)$ is a function of baseline covariates $V$ and time
  $t$. Let $g(V,t)$ be continuous in $t$ and let $e_t$ follow an
  Ornstein--Uhlenbeck process, that is,
  \[
  de_t = - \theta e_t\,dt + \sigma\,d W_t,
  \]
  where $W_t$ is the standard Brownian motion.
\item[$\bullet$]  $Y_t$ is the actual outcome process and follows the
  deterministic model (\ref{eq:dcc}):
  \[Y_t = Y^0_t +  \Psi \int^t_0 A_{s}\,ds.
  \]
\item[$\bullet$]  $A_t$ is the treatment process, taking binary values. The jump of
the $A_t$ process follows the following formula:
  \begin{eqnarray*}
   P(A_{s}\ \mbox{jumps once from }(t,t+h]|\bar{A}_t,\bar{Y}_t,\bar{Y}^0) &=& s(A_t,Y_t)h + o(h),\\
   P(A_{s}\ \mbox{jumps more than once from }(t,t+h]|\bar{A}_t,\bar{Y}_t,\bar{Y}^0) &=& o(h),
  \end{eqnarray*}
  where $\bar{A}_t$ and $\bar{Y}_t$ are the full continuous-time history
  of treatment \vspace*{1pt}and outcome up to time $t$ and $\bar{Y}^0$ is the full
  continuous-time path of potential outcome from time 0 to time $K$.
  By making $s(\cdot)$ independent of $\bar{Y}^0$, we make our model
  satisfy the continuous-time sequential randomization assumption.
\end{longlist}
In this model, the only time-dependent confounder is the outcome
process itself.

We also consider several variations of the above model (denoted as M1
below):
\begin{longlist}
\item[$\bullet$]  Model (M2) extends (M1) to the non-Markovian-$Y^0_t$ case.
  Specifically, we  consider the case where $e_t$
  in the model of $Y^0_t$ follows a non-Markovian process, namely
  an Ornstein--Uhlenbeck process in random environments, which is defined as
  the following:
  \begin{longlist}
  \item[(1)] $J_t$ is a continuous-time Markov process taking values in a
    finite set $\{1,\ldots,m\}$, which is the environment process;
  \item[(2)] we have $m>1$ sets of parameters $\theta_1, \sigma_1, \ldots, \theta_m, \sigma_m$;
  \item[(3)] $e_t$ follows an Ornstein--Uhlenbeck process with
    parameters $\theta_j,\sigma_j$, when $J_t = j$; the starting
    point of each diffusion is chosen to be simply the endpoint of the previous one.
  \end{longlist}
\item[$\bullet$] Model (M3) extends (M1) to another setting of
non-Markovian-$Y^0_t$
  process, where
  \[
  Y^0_t = g(V,t) + 0.8 e_{t-1} + 0.2 e_t.
  \]
  $e_t$ follows the same Markovian Ornstein--Uhlenbeck process as in
  M1. Every other variable is the same as in M1.
\item[$\bullet$] Model (M4) considers the case with more than one covariate. In
M4,  we keep the
  assumptions on $Y^0_t$ as in (M1) and the deterministic model of
  $Y_t$. We add one more covariate, which is generated as
  follows:
  \[
  L^-_t = 0.2 Y_t + 0.8 Y^0_{t+0.5} + 0.5 \eta_t.
  \]
  $\eta_t$ follows an  Ornstein--Uhlenbeck process independent of the
  $Y^0_t$ process.
  In this specification, the covariate $L^-_t$ contains some leading
  information about $Y^0$, but it is only ahead of $Y^0$ for 0.5
  length of a time unit. Here, we use $L^-_t$ instead of $L_t$ to denote
  that it is the covariate excluding $Y_t$. The simulation
  model for the $A_t$ process is given in Appendix \ref{appendixe}.
\end{longlist}
In all of these models, to simulate data, we use $g(V,t) = C$ (a
constant), $\Psi = 1$, a time span from 0 to 5 and a sample size of
5000. Details of other parameter specifications can be found in
Appendix \ref{appendixe}. We generate 5000 continuous paths of $Y_t$ and $A_t$ (and
$L^-_t$ in M4), from time 0 to time 5, and record $Y^*_0, A^*_0, Y^*_1,
A^*_1, \ldots, Y^*_{4}, A^*_{4}, Y^*_5$ and $\operatorname{cum}A^*_1,\ldots,\operatorname{cum}A^*_5$ (and
$L^{-*}_0,\ldots,L^{-*}_{4}$ in M4) as the observed data.

\subsection{Estimations and results under $\mathrm{M}1$}\label{s:sim:m1}
\begin{figure}

\includegraphics{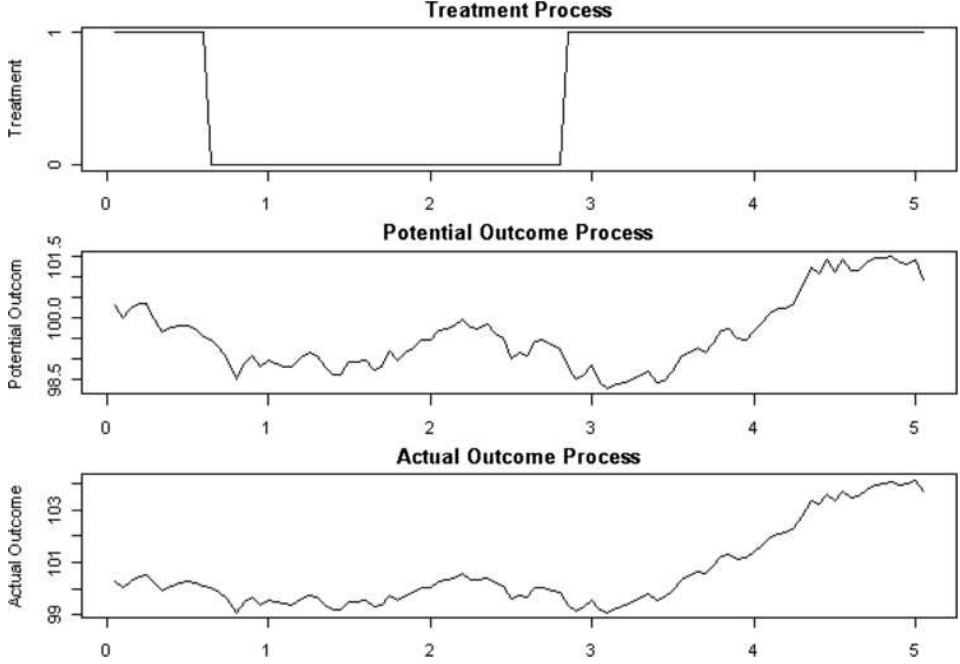}

\caption{Example of continuous-time paths under $\mathrm{M}1$.}\label{fig:process}
\end{figure}

Figure \ref{fig:process} shows a typical continuous-time path of
$Y^0_t$, $Y_t$ and $A_t$. The treatment switches around time $0.7$ and
time $2.8$.

We apply three estimating methods on data simulated from M1:
 the naive discrete-time g-estimation described in
 Section \ref{s:review:dts}, which ignores the underlying continuous-time processes;
 the modified g-estimation described in
 Section \ref{s:review:g}, which
 controls for  all the observed discrete-time history;
 and the controlling-the-future method in Section \ref{s:new:mod} of controlling for the next
period's potential outcome in addition to the discrete-time history.

For estimation, even though we know the data generating process, it is
 too complicated to use the correct model for the
propensity score,\vspace*{1.5pt} that is, the correct functional form for  $p_k(\Psi)
\equiv P(A^*_k|\bar{L}^*_k,\bar{Y}^*_k,
\bar{A}^*_{k-1},\overline{\operatorname{cum}A^*_k},Y^{0*}_{k+}(\Psi))$. Therefore, we
use the following approximations (note that we control for past
treatment and covariates as well---see comments for Theorem
\ref{thm:cond}):
\begin{longlist}
\item[(1)]  standard g-estimation ignoring continuous-time processes (naive
g-estimation)
 \[
 \operatorname{logit}(p_k) = \beta_0 + \beta_1 A^*_{k-1} + \beta_2 Y^*_{k-1} + \beta_3
  Y^*_k;
  \]
\item[(2)] g-estimation controlling for all observed history (modified
g-estimation)
  \[
  \operatorname{logit}(p_k) = \beta_0 + \beta_1 A^*_{k-1} + \beta_2 Y^*_{k-1} + \beta_3
  Y^*_k + \beta_4 \operatorname{cum}A^*_k;
  \]
\item[(3)] the controlling-the-future method, controlling for next period
potential outcomes (controlling-the-future estimation)
  \[
  \operatorname{logit}(p_k(\Psi)) = \beta_0 + \beta_1 A^*_{k-1} + \beta_2
  Y^*_{k-1} + \beta_3 Y^*_k + \beta_4 \operatorname{cum}A^*_k + \beta_5 Y^{0*}_{k+1}(\Psi).
  \]
\end{longlist}

We plug these models for the propensity scores into estimation
equations (\ref{eq:esteqs_d}), (\ref{eq:esteqdata}) and
(\ref{eq:joffeesteqs}) respectively. [Note that in equation
(\ref{eq:esteqs_d}), $Y^{0*}_k(\Psi) = Y^*_k - \Psi \sum^{k-1}_{l=0}
A^*_l$, while in the other two, $Y^{0*}_k(\Psi) = Y^*_k - \Psi
\operatorname{cum}A^*_k$.]

The first panel of Table \ref{t:m1-4} shows a summary of the estimates
of causal parameters for 1000 simulations from M1. The naive
g-estimation gives severely biased estimates. Controlling for all
observed history and controlling for additional next period potential
outcome both give us
 unbiased estimates. As discussed at the end of Section \ref{s:new:intui},
the controlling-the-future method has lower efficiency.

The last row of the first panel in Table \ref{t:m1-4} shows the coverage
rate of the 95\% confidence interval estimated from the 1000
independent simulations. Naive g-estimation has a zero coverage rate,
while the other two methods have coverage rates around 95\%.

\begin{table}
\caption{Estimated causal parameters from data generated by $\mathrm{M}1$--$\mathrm{M}4$}\label{t:m1-4}
\begin{tabular*}{\textwidth}{@{\extracolsep{\fill}}lccc@{}}
        \hline
        & \textbf{Naive g-est.} & \textbf{Mod. g-est.} & \textbf{Ctr-future est.} \\
        \hline
        \multicolumn{4}{@{}c@{}}{Simulation results from M1, \textit{true parameter $= 1$}}\\
        Mean estimate\tabnoteref[\dag]{t1}&0.7728& 1.0005  &  0.9988 \\
        S.D. of estimates\tabnoteref[\ddag]{t2}& 0.0183&0.0191&   0.0403     \\
        S.D. of the mean estimate\tabnoteref[*]{t3}& 0.0005&    0.0006 &        0.0013 \\
        Absolute bias\tabnoteref[**]{t4} &0.2272  &0.0005    &   0.0012   \\
        Coverage\tabnoteref[\diamond]{t5} & 0\phantom{0000,} & 0.946\phantom{0}& 0.956\phantom{0}
        \\ [6pt]
        \multicolumn{4}{@{}c@{}}{Simulation results from M2, \textit{true parameter $= 1$}}\\
        Mean estimate\tabnoteref[\dag]{t1}&0.7651& 1.0016  &  1.0000 \\
        S.D. of estimates\tabnoteref[\ddag]{t2}&0.0132 &  0.0158&      0.0371  \\
        S.D. of the mean estimate\tabnoteref[*]{t3}&0.0004    &0.0005 &        0.0012 \\
        Absolute bias\tabnoteref[**]{t4}&0.2349    &0.0016    &   0.0000   \\
        Coverage\tabnoteref[\diamond]{t5} &0\phantom{0000,}  & 0.953\phantom{0}& 0.950\phantom{0} \\ [6pt]
        \multicolumn{4}{@{}c@{}}{Simulation results from M3, \textit{true parameter $= 1$}}\\
        Mean estimate\tabnoteref[\dag]{t1}&0.7580& 0.9845  &   1.0026 \\
        S.D. of estimates\tabnoteref[\ddag]{t2}&0.0149&  0.0180&       0.0487 \\
        S.D. of the mean estimate\tabnoteref[*]{t3}&0.0005&   0.0006 &        0.0015 \\
        Absolute bias\tabnoteref[**]{t4} &0.2420   &0.0155   &   0.0026   \\
        Coverage\tabnoteref[\diamond]{t5} & 0\phantom{0000,} &0.855\phantom{0} & 0.956\phantom{0}\\ [6 pt]
        \multicolumn{4}{@{}c@{}}{Simulation results from M4, \textit{true parameter $= 1$}}\\
        Mean estimate\tabnoteref[\dag]{t1}&0.7816& 1.0853  &  1.0085 \\
        S.D. of estimates\tabnoteref[\ddag]{t2}&0.0201& 0.0289& 0.0806 \\
        S.D. of the mean estimate\tabnoteref[*]{t3}&0.0006&     0.0009 &  0.0025 \\
        Absolute bias\tabnoteref[**]{t4}&0.2184    &0.0853    &   0.0085   \\
        Coverage\tabnoteref[\diamond]{t5} & 0\phantom{0000} & 0.115\phantom{0}&0.948\phantom{0} \\
        \hline
\end{tabular*}
\tabnotetext[\dag]{t1}{Averaged over estimates from 1000 independent simulations of sample size 5000.}
\tabnotetext[\ddag]{t2}{Sample standard deviation of the 1000 estimates.}
\tabnotetext[*]{t3}{Sample S.D.$/\sqrt{1000}$.}
\tabnotetext[**]{t4}{Absolute value of (1-mean estimates).}
\tabnotetext[\diamond]{t5}{Coverage rate of 95\% confidence intervals for 1000 simulations.}
\end{table}

\subsection{Simulation results under $\mathrm{M}2$ and $\mathrm{M}3$}\label{s:sim:m2m3}
 The results in the second panel of Table
 \ref{t:m1-4} are
 typical for different values of parameters under M2.
 The naive g-estimation performs badly, while both of the other methods
 still work well with the data generated from M2.
 This shows that the modified g-estimation and
 the controlling-the-future method have some level of robustness to mild
 violations of the Markovian assumption.

The third part of Table \ref{t:m1-4} shows the results of simulation
from M3, where $Y^0$ violates the Markov property more substantially.
In this case, we can see that the mean of the modified g-estimates is
biased, but the mean of the controlling-the-future estimates is almost
unbiased. In the last row of the third panel, the coverage rate for the
modified g-estimation drops to
 0.855, while the controlling-the-future method still has a coverage rate of 0.956.

\subsection{Estimations and results under $\mathrm{M}4$}\label{s:sim:m4}

In M4, we create a covariate $L^-_t$ that has leading information about
$Y^0_t$. In the data simulated from M4, the observational time
sequential randomization (\ref{eq:seqrand}) no longer holds, although
the data are generated following continuous-time sequential
randomization. This simulation serves as a numerical proof of the claim
that continuous-time sequential randomization does not imply
discrete-time sequential randomization.

To show this, we consider the  following working propensity score model
at time $k=2$ and its dependence on the future potential outcome at
$m=4$:
\begin{longlist}
\item[$\bullet$] not controlling for the next period potential outcome (used in
  modified g-estimation)
  \begin{eqnarray}\label{eq:reg_nc}
    && \operatorname{logit}\bigl(P(A^*_k = 1 | \bar{A}^*_{k-}, \bar{L}^{-*}_{k},
    \overline{\operatorname{cum}A^*_{k}}, \bar{Y}^*_{k}, Y^{0*}_m)\bigr) \nonumber\\
    &&\qquad= \beta_0 + \beta_1 \operatorname{cum}A^*_k + \beta_2
    L^{-*}_{k-1} + \beta_3 L^{-*}_{k} + \beta_4 A^*_{k-1} \\
    &&\qquad\quad{}+ \beta_5 Y^*_{k} + \beta_6 Y^*_{k-1} + \beta_8 Y^{0*}_m;  \nonumber
  \end{eqnarray}
\item[$\bullet$] controlling for the next period potential outcome (used in
  controlling-the-future estimation)
  \begin{eqnarray} \label{eq:reg_c}
    && \operatorname{logit}\bigl(P(A^*_k = 1 | \bar{A}^*_{k-}, \bar{L}^{-*}_{k},
    \overline{\operatorname{cum}A^*_{k}}, \bar{Y}^*_{k}, Y^{0*}_{k+1}, Y^{0*}_m)\bigr)\nonumber \\
    &&\qquad= \beta_0 + \beta_1 \operatorname{cum}A^*_k + \beta_2
    L^{-*}_{k-1} + \beta_3 L^{-*}_{k} + \beta_4 A^*_{k-1} \\
    &&\qquad\quad{}+    \beta_5 Y^*_{k} + \beta_6 Y^*_{k-1} + \beta_7 Y^{0*}_{k+1} + \beta_8 Y^{0*}_m. \nonumber
  \end{eqnarray}
\end{longlist}
We can use the true values of  $Y^{0*}_{k+1}$ and $Y^{0*}_m$ in the
regression to test the discrete-time ignorability since we are
simulating the data. Table \ref{t:reg} shows the estimates of $\beta_7$
and $\beta_8$ in both regression models. The result shows that the
coefficient of $Y^{0*}_m$, $\beta_8$, is significant if we do not
control for the future potential outcome and is not significant if we
control for the future potential outcome. This shows that observational
time sequential randomization (\ref{eq:seqrand}) does not hold, while
the revised assumption (\ref{eq:joffe}) holds.

\begin{table}\tablewidth=255pt
\tabcolsep=0pt
\caption{Verification of observational time sequential
    randomization under $\mathrm{M}4$}\label{t:reg}
  \begin{tabular*}{255pt}{@{\extracolsep{\fill}}lcc@{}}
    \hline
    &\textbf{Reg. model (\ref{eq:reg_nc})\tabnoteref{t6}}& \textbf{Reg. model (\ref{eq:reg_c})\tabnoteref{t6}}\\
    \hline
    $\beta_7$& &0.1868\\
    $p$-value& &5.56e${-}$05\\ [3pt]
    $\beta_8$&0.0936&0.0134\\
    $p$-value&0.0006&0.691\phantom{0}\\
    \hline
   \end{tabular*}
  \tabnotetext{t6}{Simulation sample size $= 10{,}000$.}
\end{table}

The estimation results from M4 appear in the fourth panel of  Table
\ref{t:m1-4}. In applying these methods, we use the following
propensity score models separately:
\begin{longlist}
\item[(1)] g-estimation ignoring the underlying continuous-time processes
  (naive g-estimation)
  \[
  \operatorname{logit}(p_k) = \beta_0 + \beta_1 A^*_{k-1} + \beta_2 Y^*_{k-1} + \beta_3
  Y^*_k + \beta_5 L^{-*}_{k-1} + \beta_6 L^{-*}_{k};
  \]
\item[(2)] g-estimation controlling for all observed history (modified
g-estimation)
  \[
  \operatorname{logit}(p_k) = \beta_0 + \beta_1 A^*_{k-1} + \beta_2 Y^*_{k-1} + \beta_3
  Y^*_k + \beta_4 \operatorname{cum}A^*_k + \beta_5 L^{-*}_{k-1} + \beta_6 L^{-*}_{k};
  \]
\item[(3)] the controlling-the-future method controlling for next period
potential outcomes (controlling-the-future estimation)
  \begin{eqnarray*}
  \operatorname{logit}(p_k(\Psi)) &=& \beta_0 + \beta_1 A^*_{k-1} + \beta_2
  Y^*_{k-1} + \beta_3 Y^*_k + \beta_4 \operatorname{cum}A^*_k \\
  &&{}+  \beta_5L^{-*}_{k-1}+ \beta_6 L^{-*}_{k} + \beta_7 Y^{0*}_{k+1}(\Psi).
  \end{eqnarray*}
\end{longlist}

Both the naive g-estimation and the modified g-estimation give us
estimates with severe bias and they have coverage rates of 0 and 0.115,
respectively, for the 95\% confidence interval constructed from them.
It is worth noting that model 3 is misspecified, but, nevertheless,
leads to much less biased estimates, and the controlling-the-future
method has a coverage rate of 0.948.

\section{Application to the diarrhea data}\label{s:appl}

In this section, we apply the different approaches to the diarrhea
example mentioned in Section \ref{s:intro} (Example \ref{exam2}). For
illustration purposes, we ignore any informative censoring and use a
set of 224 children with complete records between ages 3 and  6 from
757 households in Bangladesh around 1998. The outcomes, $Y^*_k$, are
the heights of the children in centimeters, measured at  round $k$ of
the interviews, for $k = 1,2,3$. The treatment $A^*_k$ at the interview
$k$ is defined as $A^*_k = 1$ if  the child was sick with diarrhea
during the past two weeks of the interview and $A^*_k = 0$ otherwise.
The cumulative treatment $\operatorname{cum}A^*_k$ is the number of days that the
child suffered from diarrhea from  four months before the first
interview (July 15th, 1998) to the $k$th interview. Baseline covariates
$V$ include age in months, mother's height and whether the household
was exposed to the flood. Time-dependent covariates other than the
outcome, that is, $L^{-*}_k$, include mid-upper arm circumference,
weight for age z-score, type of toilet (open place, fixed place,
unsealed toilet, water-sealed toilet or other), garbage disposal method
(throwing away in own fixed place, throwing away in own nonfixed
place, disposing anywhere or other method), water purifying process
(filter, filter and broil, or other) and source of cooking water (from
pond or river/canal, or from tube well, ring well or supply water).

We apply naive g-estimation, modified g-estimation and the
controlling-the-future method to this data set. Since we only have
three rounds, the actual propensity score models and the estimating
equations for the three methods are as follows. Note that these
estimating equations are for illustrative purpose and may not be the
most efficient estimating equations for this data set.\vspace*{-1pt}
\begin{longlist}
\item[$\bullet$] Naive g-estimation uses the following propensity score model:\vspace*{-1pt}
\[
\operatorname{logit}\{P[A^*_{k}=1 | V, L^{-*}_{k}, Y^*_{k}]\}
=\beta_0 + \beta_V V + \beta_L L^{-*}_k + \beta_Y Y^*_k,\vspace*{-1pt}\label{eq:dihlog_n}
\]
where $k=1,2$.

The estimating equations follow the form of (\ref{eq:esteqs_d}) in
Section \ref{s:review:dts}:
\[
 \mathop{\sum_{1 \leq k<m \leq 3}}_{ 1 \leq i
\leq n} [A^*_{k,i} - P(A^*_{k,i}=1 | V_{i}, L^{-*}_{k,i}, Y^*_{k,i})]
\pmatrix{{l}Y^{0*}_{m,i}(\Psi)\vspace*{3pt}\cr  V_{i}\vspace*{3pt}\cr L^{-*}_{k,i}\vspace*{3pt}\cr
Y^*_{k,i}} = 0,\label{eq:diarrhea_n}
\]
where $Y^{0*}_{m,i}(\Psi) = Y^*_{m,i} - \Psi \sum^{m-1}_{l=1} A_l$.\vspace*{1.5pt}
\item[$\bullet$] Modified g-estimation uses this propensity score model:\vspace*{-1pt}
\begin{eqnarray*}\label{eq:dihlog_g}
&&\operatorname{logit}\{P[A^*_{k}=1 | V, L^{-*}_{k}, Y^*_{k}, \operatorname{cum}A^*_{k}]\}\\
&&\qquad=\beta_0 + \beta_V V + \beta_L L^{-*}_k + \beta_Y Y^*_k + \beta_{\operatorname{cum}A}
\operatorname{cum}A^*_k,\vspace*{-1pt}
\end{eqnarray*}
where $k=1,2$.

The estimating equations follow the form of (\ref{eq:esteqdata}) in
Section \ref{s:review:g}.
\[
 \mathop{\sum_{1\leq k < m \leq 3}}_{1 \leq i \leq
n} [A^*_{k,i} - P(A^*_{k,i}=1 | V_i, L^{-*}_{k,i}, Y^*_{k,i},
  \operatorname{cum}A^*_{k,i})]
\pmatrix{Y^{0*}_{m,i}(\Psi) \vspace*{3pt}\cr V_i \vspace*{3pt}\cr L^{-*}_{k,i}\vspace*{3pt}\cr Y^*_{k,i} \vspace*{3pt}\cr
  \operatorname{cum}A^*_{k,i}} = 0,\label{eq:diarrhea_g}
\]
where $Y^{0*}_{m,i}(\Psi) = Y^*_{m,i} - \Psi \operatorname{cum}A_m$.\vspace*{1.5pt}

\item[$\bullet$] Controlling-the-future estimation uses the following propensity
  score model:
\begin{eqnarray*}\label{eq:dihlog}
 &&\operatorname{logit}\{P[A^*_{1}=1 | V, L^{-*}_{1}, Y^*_{1},
\operatorname{cum}A^*_{1},
  Y^{0*}_{2}(\Psi)]\} \\
&&\qquad=\,\beta_0 + \beta_V V + \beta_L L^{-*}_1 + \beta_Y Y^*_1 + \beta_{\operatorname{cum}A}
\operatorname{cum}A^*_1 + \beta_{Y^0} Y^{0*}_2(\Psi).
\end{eqnarray*}

The estimating equations follow (\ref{eq:joffeesteqs}) in
Section \ref{s:new}:
\[
 \sum_{1 \leq i \leq n} [A^*_{1,i} - P(A^*_{1,i} |
V_i, L^{-*}_{1,i}, Y^*_{1,i},
  \operatorname{cum}A^*_{1,i}, Y^{0*}_{2,i}(\Psi))]
\pmatrix{Y^{0*}_{3,i}(\Psi)\vspace*{3pt}\cr V_i\vspace*{3pt}\cr L^{-*}_{1,i}\vspace*{3pt}\cr Y^*_{1,i}\vspace*{3pt}\cr
  \operatorname{cum}A^*_{1,i}\vspace*{3pt}\cr Y^{0*}_{2,i}(\Psi)} = 0,\label{eq:diarrhea}
\]
where $Y^{0*}_{3,i}(\Psi) = Y^*_{3,i} - \Psi \operatorname{cum}A_3$.
\end{longlist}

The interpretation of $\Psi$ in the last two models is that one day of
suffering from diarrhea
 reduces the height of the child by $\Psi$
centimeters. For naive g-estimation, the underlying data generating
model treats the exposure at the observational time as the constant
exposure level for the next six months, which does not make sense in
the context. It should be noted that if we apply the naive g-estimation, the estimated $\Psi$ should not
be interpreted the same way in the modified g-estimation and the controlling-the-future method.
Instead, it be interpreted as the effect of having diarrhea
at the time of visits. The effect of the child having diarrhea at any time
between the visit and the next visit six months later, but not at the time of the visit,
is not described by this $\Psi$.

The estimating equations are solved by a Newton--Raphson algorithm. The
estimated $\Psi$ and its standard deviation are reported in
Table \ref{t:dia}. Modified g-estimation estimates $\hat{\Psi} =
-0.3481$, which means that the height of the child is reduced by 0.35 cm
if the child has one day of diarrhea. Our controlling-the-future method
produces an estimate of $\hat{\Psi} = -0.0840$. Although all of the
estimates are not significant because of the small sample size, the
sign and magnitude of the estimate from the controlling-the-future
method are similar to
 what has been found in other research on diarrhea's effect on height [e.g., Moore et
al. (\citeyear{moore2001})].

\begin{table}[b]\tablewidth=220pt
\tabcolsep=0pt
   \caption{Estimation of $\Psi$ from the diarrhea data set}\label{t:dia}
  \begin{tabular*}{220pt}{@{\extracolsep{\fill}}lcc@{}}
    \hline
    \textbf{Method} & \textbf{Estimate} & \textbf{Std. err.}\\
    \hline
    Naive g-est. & $-0.3991$ & $0.2469$\\
    Modified g-est. & $-0.3481$ & $0.2832$\\
    Controlling-the-future est.& $-0.0840$ & $0.1894$\\
    \hline
  \end{tabular*}
\end{table}

In addition, we note that the standard deviation of the modified
g-estimate is higher than that of the controlling-the-future estimate.
As discussed at the end of Section \ref{s:new:intui}, if the modified
g-estimation is consistent, we would expect the controlling-the-future
estimation to have larger standard deviation. The standard deviations
in Table \ref{t:dia} provide evidence that the modified g-estimation is
not consistent.

\section{Conclusion}\label{s:concl}

In this paper, we have studied  causal inference from longitudinal data
when the underlying processes are in continuous time, but the
covariates are only observed at discrete times. We have investigated
two aspects of the problem. One is the validity of the discrete-time
g-estimation. Specifically, we investigated a modified g-estimation
that is in the spirit of standard discrete-time g-estimation, but is
modified to incorporate the information of the underlying
continuous-time treatment process, which we have  referred to as
``modified g-estimation'' throughout the paper. We have shown that an
important condition that justifies this modified g-estimation is the
finite-time sequential randomization assumption at any subset of time
points, which is strictly stronger than the continuous-time sequential
randomization. We have also shown that a Markovian assumption and the
continuous-time sequential randomization would imply the FTSR
assumption. The Markovian condition is more useful than the FTSR
assumption, in the sense that it can potentially help researchers
decide whether the application of the modified g-estimation is
appropriate.
 The other aspect is the controlling-the-future method that we propose
  to use when the condition to warrant g-estimation does not hold.
 The controlling-the-future method can produce consistent estimates
  when g-estimation is inconsistent and is less biased in other scenarios.
 In particular, we identified two important cases in which controlling
  the future is less biased, namely, when there is delayed dependence
   in the baseline potential outcome process and when there are leading
    indicators of the potential outcome process in the covariate process.

In our simulation study, we have shown the performance of the modified
g-estimation and the controlling-the-future estimation. The results
confirm our discussion in earlier sections. The simulation results also
indicate the danger of applying naive g-estimation, which is usually
severely biased and inconsistent when its underlying assumptions are
violated, as in the situations considered.

We have applied the g-estimation methods and the controlling-the-future
method to estimating the effect of diarrhea on a child's height and
estimated that its effect is negative but not significant. The real
application  also provides some evidence that the modified g-estimation
is not consistent.

All of the discussion in this paper is based on a particular form of
causal model---equation (\ref{eq:dcc}). However, all of the arguments
could apply to a class of more general rank-preserving models, with
necessary adjustments in various equations. If we assume a generic
rank-preserving model with $Y_t = f(Y^0_t, h(\bar{A}_{t-});\Psi)$,
where $\bar{A}_{t-}$ is the continuous-time path of $A$ from time $0$
to $t-$, $h$ is some functional [e.g., in our paper, $h(\bar{A}_{t-}) =
\int^t_0 A_s\,ds$] and $f$ is some strictly monotonic function with
respect to the first argument [e.g., in our paper, $f(x,y;\Psi) = x +
\Psi y$], we map $Y^*_k$ to $Y^{0*}_k =
f^{-1}(Y^*_k,h(\bar{A}_{k-});\Psi)$, where $f^{-1}$ is the inverse of
$f(x,y;\Psi)$ with respect to $x$ for any given $y$. We can then
substitute all $\operatorname{cum}A^*_k$'s in this paper by the $h(\bar{A}_{k-})$'s.
All of the discussion and formulas in the paper would remain valid
under the assumption that we observe all $h(\bar{A}_{k-})$'s, which can
be easily satisfied with detailed continuous-time records of the
treatment. It should be noted that the argument does not work  if a
time-varying covariate modifies the effect of treatment. For example,
if $Y_t = Y^0_t + \Psi \int^t_0 L^2_s A_s\,ds$, where $L_s$ is a
time-varying covariate,  observing the full continuous-time treatment
process is not enough. Some imputation for the $L_s$ process is
necessary.

The methods considered here have several limitations. These include
rank preservation, a strong assumption that the effects of treatment
are deterministic.   This assumption facilitates the interpretation of
models.  In other work on structural nested distribution and related
models [e.g., Robins (\citeyear{robins2008})], rank preservation has been shown to be
unnecessary in settings in which one is not modeling the joint
distribution of potential outcomes under different treatments. We
expect that this is also the case here,  and work justifying this more
formally is in progress. We also require that the cumulative amount of
treatment (or the full continuous-time treatment process, if using
other causal models mentioned above) between the discrete time points
when the covariates are observed is known. Work is in progress on the
more challenging case in which the  treatment process is only observed
at discrete times and the cumulative amount of treatment is measured
with error. In addition, we ignore any censoring problem requiring that
our data is complete, which might not be satisfied in reality. It will
also be interesting to study how to accommodate censored data in our
framework in future work.

\begin{appendix}
\section{\texorpdfstring{Estimating covariance matrix of estimated parameters}{Estimating
covariance matrix of estimated parameters}}\label{appendixa}
The formulas in this appendix can be used to estimate the covariance
matrix of the estimated parameters from naive g-estimation of
Section \ref{s:review:dts}, modified g-estimation of
Section \ref{s:review:g} and the controlling-the-future estimation of
Section \ref{s:new:mod}. More general results on the asymptotical
covariances can be found in van der Vaart (\citeyear{vandervaart2000}).

We write $\theta = (\Psi, \beta)$. In Sections \ref{s:review:dts} and \ref{s:review:g},
 $\beta$ is the parameter in the propensity
score model. In Section \ref{s:new:mod}, $\beta=(\beta_X,\beta_h)$ is
the parameter in the propensity score model. Let $U(\theta)$  be the
vector on the left-hand side of the estimating equations
[equation (\ref{eq:esteqs_d}) in Section \ref{s:review:dts},
equation (\ref{eq:esteqdata}) in Section \ref{s:review:g} and
equation (\ref{eq:joffeesteqs}) in Section \ref{s:new:mod},
respectively]. We also define
\[
U_{i,k,m}(\theta) \equiv \bigl(A^*_{i,k}-p_{i,k}
(\beta)\bigr)[g(Y^{0*}_{i,m}(\Psi),X^*_{i,k}),X^*_{i,k}]^T
\]
for the naive
g-estimation and the modified g-estimation, and
\[
U_{i,k,m}(\theta)\equiv \bigl(A^*_{i,k}-p_{i,k}
(\Psi;\beta_X,\beta_h)\bigr)[g(Y^{0*}_{i,m}(\Psi),X^*_i,h_i),X^*_{i,k},h_{i,k}]^T
\]
for the controlling-the-future estimation. We then h ave $U(\theta) =
\sum U_{i,k,m}$.

Let $B(\theta) = E[\frac{\partial U(\theta)}{\partial \theta}]$, which
can be estimated as
\[
 \hat{B}(\theta)=-\sum_{i,k,m} \biggl\{\frac{\partial U_{i,k,m}}{\partial
 \theta}\biggr\} \bigg|_{\theta = \hat{\theta}},
\]
where $\hat{\theta}$ is the solution from the corresponding estimating
equations,  $k<m$ in both g-estimations and $k<m-1$ in
controlling-the-future estimation. The covariance matrix of the
estimator
 $\hat{\theta}$ can then be estimated as
\[
 \operatorname{Cov}(\hat{\theta}) = \hat{B}^{-1}(\theta) \operatorname{Cov}[U(\hat{\theta})] \hat{B}^{-1}(\theta)'
\]
by the delta method, where $\operatorname{Cov}[U(\theta)]$ is estimated by
\[
\operatorname{Cov}[U(\hat{\theta})] = \sum_{i} U_{i}(\hat{\theta})U_{i}(\hat{\theta})'
\]
with $U_i = \sum_{k,m} U_{i,k,m}(\hat{\theta})$, $k<m$ in both
g-estimations and $k<m-1$ in controlling-the-future estimation.

\section{Existence of solution and identification}\label{appendixb}
The estimating equations in this paper, equation (\ref{eq:esteqs_d}) in
Section \ref{s:review:dts}, equation~(\ref{eq:esteqdata}) in
Section \ref{s:review:g} and equation (\ref{eq:joffeesteqs}) in
Section \ref{s:new:mod}, are asymptotically consistent systems of
equations by definition, if the respective underlying assumptions for
each estimating equation hold true. The existence of a solution is
guaranteed asymptotically. In addition, we have the same number of
equations as the number of parameters in each system. One would usually
expect  there to exist a solution for the estimating equations, even in
a relatively small sample.

However, the asymptotic solution may not be unique, which leads to an
identification problem. As a special case from the more general
semi-parametric theory [see Tsiatis (\citeyear{tsiatis2006})], we state the following
lemma for identification, following the notation of Appendix \ref{appendixa}.

\begin{lemma}
The parameter $\theta$ is identifiable under the model
\[
E[ U (\theta) ] = 0
\]
if both $\operatorname{Cov}[U(\theta_0)]$ and $B(\theta_0) \equiv E[\frac{\partial
U(\theta_0)}{\partial \theta}]$ are of full rank. Here, $\theta_0$ is
the value of the true parameter.
\end{lemma}

\begin{pf}
The proof is trivial. By Appendix \ref{appendixa}, the asymptotic covariance matrix
of the estimates is given by
\[
B^{-1}(\theta_0) \operatorname{Cov}[U(\theta_0)] B^{-1}(\theta_0)',
\]
which will be finite and of full rank when the conditions in the lemma hold
true.
\end{pf}\vspace*{-15pt}

\section{Proof that FTSR implies CTSR}\label{appendixc}
We assume that $Z_t$ is a c\`{a}dl\`{a}g process, and everything we
discuss is in an a.s. sense.

We first define
\begin{eqnarray*}
\mathcal{H}_{t-} & \equiv &\sigma(\bar{Z}_{t-}), \\
\mathcal{F}_{t-,t+} & \equiv &\sigma(\bar{Z}_{t-},\underline{Y}^0_{t+}).
\end{eqnarray*}

Recall that $N_t$ counts the number of jumps in $A_t$ up to time $t$.
We assume that a continuous version of the $\mathcal{F}_{t-,t+}$
intensity process of $N_t$ exists, which we denote by $\eta_t$. If we
define
\[
r_t(\delta) = (1 - A_{t-}) A_{t+\delta} + A_{t-} (1-A_{t+\delta}).
\]

Then, under certain regularity conditions [see Chapter 2 of Andersen et
al. (\citeyear{andersen1992})], for every $t$,
\[
\eta_t = \lim_{\delta \downarrow 0}
\frac{E[r_t(\delta)|\mathcal{F}_{t-,t+}]}{\delta}\qquad\mbox{a.s.}
\]

For Theorem \ref{thm:main}, we need to show that $\eta_t$ is also
$\mathcal{H}_{t-}$-measurable. This is because if this is true, then
\begin{eqnarray*}
&&E\biggl[N_{t_0+s} - \int^{t_0+s}_0 \eta_t\,dt \Big| \mathcal{H}_{t_0}\biggr] \\
&&\qquad=  E\biggl[\biggl(N_{t_0+s} - \int^{t_0+s}_0 \eta_t\,dt\biggr) - \biggl(N_{t_0} - \int^{t_0}_0 \eta_t\,dt\biggr)
\Big| \mathcal{H}_{t_0}\biggr] \\
&&\qquad\quad{}+ E\biggl[\biggl(N_{t_0} - \int^{t_0}_0 \eta_t\,dt\biggr) \Big| \mathcal{H}_{t_0}\biggr] \\
 &&\qquad=  E\biggl\{E\biggl[\biggl(N_{t_0+s} - \int^{t_0+s}_0 \eta_t\,dt\biggr) - \biggl(N_{t_0} - \int^{t_0}_0 \eta_t\,dt\biggr)
 \Bigl| \mathcal{F}_{t_0-,t_0+}\biggr] \Bigr|\mathcal{H}_{t_0}\biggr\} \\
 &&\qquad\quad{}  + N_{t_0} - \int^{t_0}_0 \eta_t\,dt \\
&&\qquad = 0 + N_{t_0} - \int^{t_0}_0 \eta_t\,dt \\
 &&\qquad= N_{t_0} - \int^{t_0}_0 \eta_t\,dt.
\end{eqnarray*}

The second equality follows because of properties of conditional
expectation and the assumption that $\eta_t$ is
$\mathcal{H}_{t-}$-measurable. The third equality holds because
$\eta_t$ is an $\mathcal{F}_{t-,t+}$ intensity process of $N_t$. The
last equality shows that $\eta_t$ is also a  $\mathcal{H}_{t-}$
intensity process of $N_t$, which agrees with the definition of CTSR.

Before proving the main result, we assume the following regularity
conditions.
\begin{longlist}
\item[1.] As stated before, we assume that $\eta_t$ is continuous.  We
further assume that $\eta_t$ is positive, and bounded from below and
above by constants that do not depend on $t$. We also assume that
$\frac{E[r_t(\delta)|\mathcal{F}_{t-,t+}]}{\delta}$ is bounded by a
constant for every $t$ within a interval of $(0, \delta_0]$.
\item[2.] We
assume that for any finite sequence of time points, $t_1 \leq t_2 \leq
t_3 \leq \cdots \leq t_n$, the density $f(Z_{t_1} = z_1, Z_{t_2} = z_2,
\ldots, Z_{t_n} = z_n)$ is well defined and  locally uniformly bounded,
that is, there exists a constant $D$ and a rectangle $B \equiv [t_1 -
\delta_1, t_1 + \delta_1] \times [t_2 - \delta_2, t_2 + \delta_2]
\times \cdots \times [t_n - \delta_n, t_n + \delta_n]$ such that for
any $(t'_1, t'_2, \ldots, t'_n)^T \in B$ and
     any possible value of $(z_1, z_2, \ldots, z_n)^T$,
    \[
    f(Z_{t'_1} = z_1, Z_{t'_2} = z_2, \ldots, Z_{t'_n} = z_n) \leq D.
    \]
    For any conditional expectation involving finite sequence of time points, we choose the version that is defined by the joint density.
\item[3.] Given any finite sequence of time points, $t_1 \leq t_2 \leq t_3
\leq \cdots \leq t_n$ and any possible value of $(z_1, z_2, \ldots,
z_n)^T$, we assume that the following convergence is uniform in a
closed neighborhood of $\tilde{t} \equiv (t_1 , t_2 , t_3 ,\ldots ,
t_n)$:
    \begin{eqnarray*}
    &&f(Z_{t'_1} = z_1, Z_{t'_2} = z_2, \ldots, Z_{t'_n} = z_n) \\
    &&\qquad=  \lim_{\Delta \downarrow 0} \frac{P(Z_{t'_1}
     \in [z_1, z_1+\Delta_1], Z_{t'_2} \in [z_2, z_2 + \Delta_2],
      \ldots, Z_{t'_n} \in [z_n, z_n+\Delta_n])}{\Delta_1 \times \Delta_2 \times \cdots \times \Delta_n},
    \end{eqnarray*}
    where $(t'_1, t'_2, \ldots, t'_n)^T$ is in a neighborhood of $\tilde{t}$.
\item[4.] Given any finite sequence of time points, $t_1 \leq t_2 \leq t_3
\leq \cdots \leq t_i \leq \cdots \leq t_n$ and any possible value of
$(z_1, z_2, \ldots, z_n)^T$, we define
    \[
    f(\delta) = \frac{P(A_{t_i+\delta} \neq A_{t_i}|Z_{t_1} = z_1, Z_{t_2} = z_2, \ldots, Z_{t_n} = z_n)}{\delta}.
    \]
    We assume that $\lim_{\delta \downarrow 0} f(\delta)$ exists and is positive and finite. We also assume that $f(\delta)$ is finite and is right-continuous in $\delta$, and the continuity is uniform with respect
    to $(\delta, t_i)$ in $[0,\delta_0] \times B(t_i)$, where $B(t_i)$ is a closed neighborhood of $t_i$. Further, we assume that the above assumption is true if any of the $Z$ in $f$ is in its left-limit value rather than the concurrent value.
\end{longlist}

\begin{remark}\label{rem:regularity3}
 The third regularity condition is needed when
we want to prove convergence in density. For example, consider that
when $\delta \downarrow 0$, we have $Z_{t_2+\delta} \rightarrow
Z_{t_2}$. We can then see that
\begin{eqnarray*}
&&\hspace*{-4.8pt}\lim_{\delta \downarrow 0} f(Z_{t_1} = z_1, Z_{t_2+\delta} = z_2, Z_{t_3} = z_3) \\
&&\!\hspace*{-4.8pt}\qquad= \lim_{\delta \downarrow 0}\mathop{\mathop{\lim_{\Delta_1 \downarrow 0}}_{\Delta_2 \downarrow 0}}_{\Delta_3 \downarrow 0} \frac{P(Z_{t_1} \in [z_1, z_1 + \Delta_1], Z_{t_2+\delta} \in [z_2, z_2 + \Delta_2], Z_{t_3} \in [z_3, z_3+\Delta_3])}{\Delta_1 \Delta_2 \Delta_3} \\
&&\!\hspace*{-4.8pt}\qquad=\mathop{\mathop{\lim_{\Delta_1 \downarrow 0}}_{\Delta_2 \downarrow 0}}_{\Delta_3 \downarrow 0}\lim_{\delta \downarrow 0} \frac{P(Z_{t_1} \in [z_1, z_1 + \Delta_1], Z_{t_2+\delta} \in [z_2, z_2 + \Delta_2], Z_{t_3} \in [z_3, z_3+\Delta_3])}{\Delta_1 \Delta_2 \Delta_3} \\
&&\!\hspace*{-4.8pt}\qquad=\mathop{\mathop{\lim_{\Delta_1 \downarrow 0}}_{ \Delta_2 \downarrow 0}}_{\Delta_3 \downarrow 0}\frac{P(Z_{t_1} \in [z_1, z_1 + \Delta_1], Z_{t_2} \in [z_2, z_2 + \Delta_2], Z_{t_3} \in [z_3, z_3+\Delta_3])}{\Delta_1 \Delta_2 \Delta_3} \\
&&\!\hspace*{-4.8pt}\qquad= f(Z_{t_1} = z_1, Z_{t_2} = z_2, Z_{t_3} = z_3).
\end{eqnarray*}
The  interchanging of limits in the second equality is valid because of
the third regularity condition. The third equality follows from the
fact that probabilities are expectations of indicator functions and
that the dominated convergence theorem applies.
\end{remark}

We introduce the following lemma for technical convenience.
\begin{lemma}\label{lem:unlimit}
 If the \textit{c\`{a}dl\`{a}g} process $Z_t$
follows the \textit{finite-time
  sequential randomization} as defined in Definition \ref{def:fti},
then the following version of \textit{FTSR} is also true:
\begin{eqnarray} \label{eq:fti_m}
  &&P(A_{t_n} | \bar{L}_{t_{n-1}},L_{t_n-}, \bar{A}_{t_{n-1}},
    \bar{Y}^0_{t_{n-1}}, Y^0_{t_n-},
    \underline{Y}^0_{t_n+}) \nonumber\\ [-8pt]\\ [-8pt]
     &&\qquad=P(A_{t_n} | \bar{L}_{t_{n-1}}, L_{t_n-},
    \bar{A}_{t_{n-1}}, \bar{Y}^0_{t_{n-1}}, Y^0_{t_n-}),\nonumber
\end{eqnarray}
where $\bar{L}_{t_{n-1}} = (L_{t_1},L_{t_2},\ldots,L_{t_{n -1} })$,
$\bar{A}_{t_{n-1}} = (A_{t_1},A_{t_2},\ldots,A_{t_{n-1}})$,
$\bar{Y}^0_{t_{n-1}} = (Y^0_{t_1},\break Y^0_{t_2},\ldots,Y^0_{t_{n -1} })$
and $\underline{Y}^0_{t_n+} =
(Y^0_{t_{n+1}},Y^0_{t_{n+2}},\ldots,Y^0_{t_{n+l}})$.
\end{lemma}

\begin{remark}
The difference between (\ref{eq:fti_m}) and the original definition of
FTSR is that in (\ref{eq:fti_m}), most $L$'s and $Y^0$'s are stated in
  their concurrent values, while in Definition \ref{def:fti}, they are
  all stated in their left limits. Lemma \ref{lem:unlimit} is only for
  technical convenience.
\end{remark}

\begin{pf*}{Proof of Lemma \ref{lem:unlimit}}
The result follows directly from the definition of a c\`{a}dl\`{a}g
process.
\end{pf*}

We now consider a discrete-time property.
\begin{lemma}\label{lem:nfti}
 Suppose FTSR holds true. If we define
\begin{eqnarray*}
\mathcal{F} & =& \sigma(Z_{t_1},\ldots,Z_{t_{n-1}}, Z_{t-},
Y^0_{t_{n+1}},\ldots, Y^0_{t_{n+l}}), \\
\mathcal{H} & = &\sigma(Z_{t_1},\ldots,Z_{t_{n-1}}, Z_{t-}),
\end{eqnarray*}
then we have for every $t$ that
\begin{equation}\label{eq:nfti}
\lim_{\delta \downarrow 0} \frac{E[r_t(\delta) |
    \mathcal{F}]}{\delta} = \lim_{\delta \downarrow 0} \frac{E[r_t(\delta) |
    \mathcal{H}]}{\delta}\qquad a.s.
\end{equation}
\end{lemma}

\begin{pf}
First, we note that the limits on both sides of equation
(\ref{eq:nfti}) exist and are finite. This fact follows from the
regularity condition 1. Take $\lim_{\delta \downarrow 0}
\frac{E[r_t(\delta) |    \mathcal{F}]}{\delta}$, for example:
\begin{eqnarray*}
\lim_{\delta \downarrow 0} \frac{E[r_t(\delta) |  \mathcal{F}]}{\delta}&=&\lim_{\delta \downarrow 0} \frac{E[E[r_t(\delta)|\sigma(\bar{Z}_{t-},\underline{Y}^0_{t})] |    \mathcal{F}]}{\delta} \\
&=&E\biggl[\lim_{\delta \downarrow 0} \frac{E[r_t(\delta)|\sigma(\bar{Z}_{t-},\underline{Y}^0_{t})] }{\delta}\Big|\mathcal{F}\biggr]\\
&=&E[\eta_t | \mathcal{F}].
\end{eqnarray*}
The interchange of limit and expectation is guaranteed by the
assumption in regularity condition 1 that
$\frac{E[r_t(\delta)|\sigma(\bar{Z}_{t-},\underline{Y}^0_{t})]
}{\delta}$ is bounded. The existence is then guaranteed by the
dominated convergence theorem and $E[\eta_t | \mathcal{F}]$ is
obviously finite.

Given equation (\ref{eq:fti}) and Lemma \ref{lem:unlimit}, we always
have
\begin{equation}
\label{eq:fti2} E[I_{A_t \neq A_{t_{n}}} | \bar{L}_{t-},
  \bar{A}_{t_{n}},\bar{Y}^0_{t-},\underline{Y}^0_{t+}] = E[I_{A_t \neq
  A_{t_{n}}} | \bar{L}_{t-},
  \bar{A}_{t_{n}},\bar{Y}^0_{t-}]\qquad\mbox{a.s.},
\end{equation}
where $\bar{L}_{t-} = (L_{t_1},L_{t_2},\ldots,L_{t_{n-1}},L_{t_n},
L_{t -})^T$, $\bar{A}_{t_{n}} = (A_{t_1},A_{t_2},\ldots,A_{t_{n}})^T$,
$\bar{Y}^0_{t-} =(Y^0_{t_1},Y^0_{t_2},\ldots,Y^0_{t_{n}},Y^0_{t-})^T$
and
$\underline{Y}^0_{t+}=(Y^0_{t_{n+1}},Y^0_{t_{n+2}},\ldots,Y^0_{t_{n+l}})^T$.\vspace*{1.5pt}

In the regularity conditions, since we assumed the existence of joint
density, the usual definition of conditional probability is a version
of the conditional expectation defined using $\sigma$-fields. In our
case, we have
\begin{eqnarray*}
&&\lim_{\delta \downarrow 0} \frac{E[r_t(\delta) |
    \mathcal{F}]}{\delta} \\
&&\qquad= \lim_{\delta \downarrow 0}
    \frac{P(A_{t+\delta} \neq A_{t-} |  Z_{t_1},\ldots,Z_{t_{n-1}}, Z_{t-},
Y^0_{t_{n+1}},\ldots, Y^0_{t_{n+l}})}{\delta}\\
&&\qquad= \lim_{\delta \downarrow 0} \lim_{t_n \uparrow t-}
 \frac{P(A_{t+\delta} \neq A_{t_n} |  Z_{t_1},\ldots,Z_{t_{n-1}},
    Z_{t_n}, L_{t-}, Y^0_{t-},
Y^0_{t_{n+1}},\ldots, Y^0_{t_{n+l}})}{\delta + (t-t_n)}\\
&&\qquad=  \lim_{t_n \uparrow t-}  \lim_{\delta \downarrow 0}
 \frac{P(A_{t+\delta} \neq A_{t_n} |  Z_{t_1},\ldots,Z_{t_{n-1}},
    Z_{t_n}, L_{t-}, Y^0_{t-},
Y^0_{t_{n+1}},\ldots, Y^0_{t_{n+l}})}{\delta + (t-t_n)}\\
&&\qquad=  \lim_{t_n \uparrow t-}
 \frac{P(A_{t} \neq A_{t_n} |  Z_{t_1},\ldots,Z_{t_{n-1}},
    Z_{t_n}, L_{t-}, Y^0_{t-},
Y^0_{t_{n+1}},\ldots, Y^0_{t_{n+l}})}{ t-t_n}\\
&&\qquad=  \lim_{t_n \uparrow t-} \frac{E[I_{A_t \neq A_{t_{n}}} |
\bar{L}_{t-},
  \bar{A}_{t_{n}},\bar{Y}^0_{t-},\underline{Y}^0_{t+}]}{t - t_n}.
\end{eqnarray*}
The second equality is guaranteed by the third regularity condition. By
Remark \ref{rem:regularity3}, we can show that the conditional density
in the third line converges to the second line as $A_{t_n}$ and
$Z_{t_n}$ converges to $A_{t-}$ and $Z_{t-}$. The $(t-t_n)$ term in the
denominator is not needed for the second equality, but is crucial for
the interchangeability of limits in the third equality. The
interchangeability of limits is guaranteed by the fourth regularity
condition. By the fourth regularity condition,  the following limit
\begin{eqnarray*}
&&\lim_{\delta \downarrow 0}
 \frac{P(A_{t+\delta} \neq A_{t_n} |  Z_{t_1},\ldots,Z_{t_{n-1}},
    Z_{t_n}, Z_{t-},
Z_{t_{n+1}},\ldots, Z_{t_{n+l}})}{\delta + (t-t_n)} \\
&&\qquad=\frac{P(A_{t} \neq A_{t_n} |  Z_{t_1},\ldots,Z_{t_{n-1}},
    Z_{t_n}, Z_{t-},
Z_{t_{n+1}},\ldots, Z_{t_{n+l}})}{ t-t_n}
\end{eqnarray*}
is uniform in $t_n$.

If we integrate out some extra variables, we can get that
\begin{eqnarray*}
&&\lim_{\delta \downarrow 0}
 \frac{P(A_{t+\delta} \neq A_{t_n} |  Z_{t_1},\ldots,Z_{t_{n-1}},
    Z_{t_n}, L_{t-}, Y^0_{t-},
Y^0_{t_{n+1}},\ldots, Y^0_{t_{n+l}})}{\delta + (t-t_n)} \\
&&\qquad=\frac{P(A_{t} \neq A_{t_n} |  Z_{t_1},\ldots,Z_{t_{n-1}},
    Z_{t_n}, L_{t-}, Y^0_{t-},
Y^0_{t_{n+1}},\ldots, Y^0_{t_{n+l}})}{ t-t_n}
\end{eqnarray*}
is uniform in $t_n$.

Therefore, we can interchange the limits in the third equality.

Similarly, we can prove that
\[
\lim_{\delta \downarrow 0} \frac{E[r_t(\delta) |
    \mathcal{H}]}{\delta} = \lim_{t_n \uparrow t-}
\frac{E[I_{A_t \neq A_{t_{n}}} | \bar{L}_{t-},
  \bar{A}_{t_{n}},\bar{Y}^0_{t-}]} {t - t_n}.
\]

Therefore, we have
\begin{eqnarray*}
\lim_{\delta \downarrow 0} \frac{E[r_t(\delta) |
    \mathcal{F}]}{\delta} &=& \lim_{t_n \uparrow t-}
\frac{E[I_{A_t \neq A_{t_{n}}} | \bar{L}_{t-},
  \bar{A}_{t_{n}},\bar{Y}^0_{t-},\underline{Y}^0_{t+}]} {t - t_n} \\
&=& \lim_{t_n \uparrow t-} \frac{E[I_{A_t \neq A_{t_{n}}} |
\bar{L}_{t-},
  \bar{A}_{t_{n}},\bar{Y}^0_{t-}]} {t - t_n} \\
&=& \lim_{\delta \downarrow 0} \frac{E[r_t(\delta) |
    \mathcal{H}]}{\delta}.
\end{eqnarray*}
The second equality comes from (\ref{eq:fti2}).
\end{pf}

We now prove the final key lemma.
\begin{lemma}\label{lem:main}
 Given FTSR, $\eta_t$ is $\mathcal{H}_{t-}$-measurable.
\end{lemma}

\begin{pf}
We prove the result by using the definition of a measurable function
with respect to a $\sigma$-field.

For any $a \in \mathcal{R}$, consider the following set:
\[
B \equiv \{\omega\dvtx \eta_t \leq a\}.
\]

Since $\eta_t$ is measurable with respect to $\mathcal{F}_{t-,t+}$, $B
\in \mathcal{F}_{t-,t+}$.

By Lemma 25.9 of Rogers and Williams (\citeyear{rogers1994}), $B$ is a $\sigma$-cylinder
and it can be decided by variables from countably many time points.
Suppose the collection of these countably many time points is $S$. $S =
S_1 \cup S_2$, where $t_{1,i} < t$ for $t_{1,i} \in S_1$ and
$t_{2,j}
> t$ for $t_{2,j} \in S_2$.

Let $\mathcal{F}_{S}$ denote the $\sigma$-field generated by
$(Z_{t_{1,i}}, i \in \mathcal{N}; Z_{t-}; Y^0_{t_{2,j}}, j \in
\mathcal{N})$. We have augmented the $\sigma$-field generated by
variables from $S$ with $Z_{t-}$.

Next, define the following series of $\sigma$-fields:
\begin{eqnarray*}
\mathcal{F}_1 & \equiv& \sigma(Z_{t_{1,1}},Z_{t-},Y^0_{t_{2,1}}), \\
\mathcal{F}_2 & \equiv& \sigma(\mathcal{F}_1, Z_{t_{1,2}},Y^0_{t_{2,1}}), \\
&\cdots& \\
\mathcal{F}_{\infty} & \equiv& \mathcal{F}_{S}.
\end{eqnarray*}

Considering the following sets:
\begin{eqnarray*}
B_1 &\equiv& \{\omega\dvtx E[\eta_t | \mathcal{F}_1] \leq a\}, \\
B_2 &\equiv& \{\omega\dvtx E[\eta_t | \mathcal{F}_2] \leq a\}, \\
&\cdots& \\
B_S \equiv B_{\infty} &\equiv& \{\omega\dvtx E[\eta_t | \mathcal{F}_{S}]\leq a\}.
\end{eqnarray*}

We have $B_k \in \mathcal{F}_k$.

It is easy to see that
\[
B_1 \supset B_2 \supset \cdots \supset B_{S}
\]
because
\[
 E[ E[\eta_t | \mathcal{F}_k] | \mathcal{F}_{k-1} ] = E[\eta_t |
  \mathcal{F}_{k-1}]
\]
and taking conditional expectation preserves the direction of
inequality.

Also, with the above definitions, $\mathcal{F}_k \uparrow
\mathcal{F}_S$. Therefore, by Theorem 5.7 from Durrett [(\citeyear{durrett2005}), Chapter
4], we know that
\[
E[\eta_t | \mathcal{F}_k] \rightarrow E[\eta_t | \mathcal{F}_S]
\qquad\mbox{a.s.}
 \]
 It is then easy to see that $I_{B_1} \rightarrow
I_{B_S}$ a.s. and that
\[
B_{S} = \bigcap^{\infty}_{i=1} B_i
\]
with difference up to a null set.

We now claim that
\begin{equation}\label{eq:sets}
 B_{S} = B
\end{equation}
with difference up to a null set.

Obviously, $B \subset B_{S}$. Suppose that $P(B_{S} - B) > 0$. Since
$B_{S} - B \in \mathcal{F}_{S}$, we have
\[
\int_{B_{S}-B} \eta_t P(d\omega) = \int_{B_{S}-B} E[\eta_t |
  \mathcal{F}_S] P(d\omega).
\]

Then
\[
\mathit{LHS} > a P(B_{S} - B)
\]
and
\[
\mathit{RHS} \leq a P(B_S - B).
\]

This is a contradiction.

Therefore, $B = \bigcap^{\infty}_{i=1} B_i$ with difference up to a
null set.

Next, we define
\begin{eqnarray*}
\mathcal{H}_1 &=& \sigma(Z_{t_{1,1}}, Z_{t-}), \\
\mathcal{H}_2 &=& \sigma(\mathcal{H}_1, Z_{t_{1,2}}), \\
&\cdots&
\end{eqnarray*}

Given FTSR, by Lemma \ref{lem:nfti}, we have
\[
E[\eta_t | \mathcal{F}_k] = E[\eta_t | \mathcal{H}_k].
\]

Therefore, every $B_k \in \mathcal{H}_k$ and thus $B_k \in
\mathcal{H}_{t-}$.

Since $B = \bigcap^{\infty}_{i=1} B_i$, $B \in \mathcal{H}_{t-}$ as
well. By the definition of a measurable function, $\eta_t$ is
measurable with respect to $\mathcal{H}_{t-}$.
\end{pf}

Combining all of the results in this appendix, we have proven
Theorem \ref{thm:main}.

\section{\texorpdfstring{Proof of Theorem \protect\lowercase{\ref{thm:cond}}}{Proof of Theorem 5}}
\label{appendixd}

Let $\mathcal{G}_t = \sigma(Y^0_{t-}, L_{t-}, A_{t-})$. Recall the
definition of $r_t(\delta) = (1 - A_{t-})A_{t+\delta} +
A_{t-}(1-A_{t+\delta})$ and that $Z_t = (Y^0_t, L_t, A_t)^T$. By the
Markovian property and the c\'{a}dl\'{a}g property, it is easy to show
that
\[
E[r_t(\delta) | \sigma(\bar{Z}_{t-})] = E[r_t(\delta) | \mathcal{G}_{t}]
\]
and that
\[
E[r_t(\delta) | \sigma(\bar{Z}_{t-}, Y^0_{t+s})] = E[r_t(\delta) | \sigma(\mathcal{G}_{t}, Y^0_{t+s})].
\]
Note that, without loss of generality, we only consider $Y^0_{t+s}$ in
the proof, rather than $\underline{Y}^0_{t+}$.

Therefore, we have a reduced form of \textit{continuous-time sequential
randomization}:
\begin{eqnarray*}
\lim_{\delta \downarrow 0} \frac{E[r_t(\delta) | \sigma(\mathcal{G}_{t}, Y^0_{t+s})]}{\delta} &=& \lim_{\delta \downarrow 0} \frac{E[r_t(\delta) | \sigma(\bar{Z}_{t-}, Y^0_{t+s})]}{\delta} \\
&=& \lim_{\delta \downarrow 0} \frac{E[r_t(\delta) | \sigma(\bar{Z}_{t-})]}{\delta} \\
&=&\lim_{\delta \downarrow 0} \frac{E[r_t(\delta) |
\mathcal{G}_{t}]}{\delta}.
\end{eqnarray*}

First, we note that if we can prove
\begin{eqnarray}\label{eq:pfmarkov1}
&&f(Y^0_{t+s}, A_{t-} | Y^0_{t-}, L_{t-}) P(A_t |Y^0_{t-},
L_{t-})\nonumber\\ [-8pt]\\ [-8pt]
&&\qquad =f(Y^0_{t+s}, A_{t} | Y^0_{t-}, L_{t-}) P(A_{t-} |Y^0_{t-},
L_{t-}),\nonumber
\end{eqnarray}
then we can conclude (\ref{eq:markovcon}). The reason is as follows:
assuming  (\ref{eq:pfmarkov1}) to be true, we integrate $A_{t-}$ out on
both sides of the equation. We will get
\[
f(Y^0_{t+s}| Y^0_{t-}, L_{t-}) P(A_t |Y^0_{t-}, L_{t-}) = f(Y^0_{t+s}, A_{t} | Y^0_{t-}, L_{t-}).
\]
Dividing the above equation by $f(Y^0_{t+s}| Y^0_{t-}, L_{t-})$, we
obtain (\ref{eq:markovcon}).

Consider
\[
g(\delta_1, \delta_2) \equiv f(Y^0_{t+s}|A_{t+\delta_1}=a_1, A_{t-\delta_2}=a_2, Y^0_{t-}, L_{t-}),
\]
where $\delta_1 > 0$ and $\delta_2 > 0$.

We observe that
\begin{eqnarray*}
&&\hspace*{-4.6pt}\lim_{\delta_1 \downarrow 0} \lim_{\delta_2 \downarrow 0} g(\delta_1, \delta_2) \\
&&\qquad\hspace*{-4.6pt}=\lim_{\delta_1 \downarrow 0} f(Y^0_{t+s} | A_{t+\delta_1}=a_1 , A_{t-}=a_2,  Y^0_{t-}, L_{t-}) \\
&&\qquad\hspace*{-4.6pt}= \lim_{\delta_1 \downarrow 0}\frac{f(Y^0_{t+s}, A_{t+\delta_1}=a_1 | A_{t-}=a_2,  Y^0_{t-}, L_{t-})}{P(A_{t+\delta_1} | A_{t-}=a_2,  Y^0_{t-}, L_{t-})} \\
&&\qquad\hspace*{-4.6pt}= f(Y^0_{t+s}| A_{t-}=a_2,  Y^0_{t-}, L_{t-})\\
&&\qquad\hspace*{-4.6pt}\quad{}\times\lim_{\delta_1 \downarrow 0}\frac{P(A_{t+\delta_1}=a_1 | Y^0_{t+s}, A_{t-}=a_2,  Y^0_{t-}, L_{t-})}
{P(A_{t+\delta_1}=a_1 | A_{t-}=a_2,  Y^0_{t-}, L_{t-})} \\
&&\qquad\hspace*{-4.6pt}= \cases{
\displaystyle f(Y^0_{t+s}| A_{t-}=a_2,  Y^0_{t-}, L_{t-})\cr
\qquad\displaystyle {}\times \lim_{\delta_1 \downarrow 0}\frac{1-P(A_{t+\delta_1}
\neq A_{t-} | Y^0_{t+s}, A_{t-}=a_2,  Y^0_{t-}, L_{t-})}{1-P(A_{t+\delta_1}\neq A_{t-} | A_{t-}=a_2,  Y^0_{t-},
L_{t-})},\cr
\qquad \mbox{if } a_1 = a_2,\cr
\displaystyle f(Y^0_{t+s}| A_{t-}=a_2,  Y^0_{t-}, L_{t-})\cr
\qquad\displaystyle{} \times \lim_{\delta_1 \downarrow 0}\frac{P(A_{t+\delta_1}
\neq A_{t-} | Y^0_{t+s}, A_{t-}=a_2,  Y^0_{t-}, L_{t-})/\delta_1}
{P(A_{t+\delta_1}\neq A_{t-} | A_{t-}=a_2,  Y^0_{t-},
L_{t-})/\delta_1},\cr
\qquad \mbox{if } a_1 \neq a_2.}\\
&&\qquad\hspace*{-4.6pt}= f(Y^0_{t+s}| A_{t-}=a_2,  Y^0_{t-}, L_{t-}).
\end{eqnarray*}
Here, the validity of taking the limit inside the density is guaranteed
by the third regularity condition, and the last equality follows
because of the continuous-time sequential randomization assumption.

We also observe that
\begin{eqnarray*}
\lim_{\delta_2 \downarrow 0} \lim_{\delta_1 \downarrow 0} g(\delta_1, \delta_2)
&=&\lim_{\delta_2 \downarrow 0} f(Y^0_{t+s} | A_{t-\delta_2}, A_{t},  Y^0_{t-}, L_{t-}) \\
&=&\lim_{\delta_2 \downarrow 0} f(Y^0_{t+s} | A_{t},  Y^0_{t-}, L_{t-})=f(Y^0_{t+s} | A_{t},  Y^0_{t-}, L_{t-}).
\end{eqnarray*}

The second equality uses the Markov property.

If we can interchange the limits, then we have
\[
f(Y^0_{t+s}| A_{t-},
Y^0_{t-}, L_{t-}) = f(Y^0_{t+s} | A_{t},  Y^0_{t-}, L_{t-}).
\]
Equation
(\ref{eq:pfmarkov1}) follows from the definition of conditional
density.

We now establish the fact that
\[
\lim_{\delta_2 \downarrow 0} \lim_{\delta_1 \downarrow 0} g(\delta_1, \delta_2)
= \lim_{\delta_1 \downarrow 0} \lim_{\delta_2 \downarrow 0} g(\delta_1, \delta_2)
\]
by showing that $\lim_{\delta_1 \downarrow 0} g(\delta_1, \delta_2)$ is
uniform in $\delta_2$.

If we define $g_1(\delta_2) = \lim_{\delta_1 \downarrow 0} g(\delta_1,
\delta_2)$, then
\begin{eqnarray*}
|g(\delta_1, \delta_2) - g_1(\delta_2)|&=&\biggl|\frac{f(Y^0_{t+s},
A_{t+\delta_1}=a_1| A_{t-\delta_2}=a_2, Y^0_{t-}, L_{t-})}{P(A_{t+\delta_1}=a_1| A_{t-\delta_2}=a_2, Y^0_{t-}, L_{t-})} \\
&&\hspace*{4pt}{}- \frac{f(Y^0_{t+s},A_{t}=a_1| A_{t-\delta_2}=a_2, Y^0_{t-}, L_{t-})}{P(A_{t}=a_1| A_{t-\delta_2}=a_2, Y^0_{t-}, L_{t-})}\biggr| \\
&=& f(Y^0_{t+s}|A_{t-\delta_2}=a_2, Y^0_{t-}, L_{t-}) \\
&&{}\times \biggl|\frac{P(A_{t+\delta_1}=a_1| A_{t-\delta_2}=a_2,
Y^0_{t-}, L_{t-},Y^0_{t+s})}{P(A_{t+\delta_1}=a_1| A_{t-\delta_2}=a_2, Y^0_{t-}, L_{t-})} \\
&&\hspace*{17pt}{}- \frac{P(A_{t}=a_1| A_{t-\delta_2}=a_2, Y^0_{t-},
L_{t-},Y^0_{t+s})}{P(A_{t}=a_1| A_{t-\delta_2}=a_2, Y^0_{t-},
L_{t-})}\biggr|.
\end{eqnarray*}

Consider the ratio $\frac{P(A_{t+\delta_1}=a_1| A_{t-\delta_2}=a_2,
Y^0_{t-}, L_{t-},Y^0_{t+s})}{P(A_{t+\delta_1}=a_1| A_{t-\delta_2}=a_2,
Y^0_{t-}, L_{t-})}$. We claim that it converges to $\frac{P(A_{t}=a_1|
A_{t-\delta_2}=a_2, Y^0_{t-}, L_{t-},Y^0_{t+s})}{P(A_{t}=a_1|
A_{t-\delta_2}=a_2, Y^0_{t-}, L_{t-})}$ uniformly in $\delta_2$.

If $a_1 = a_2$, then the density $P(A_{t+\delta_1}=a_1|
A_{t-\delta_2}=a_2, Y^0_{t-}, L_{t-})$ is bounded from below by a
positive number. By the fourth regularity condition,
\begin{eqnarray*}
&&P(A_{t+\delta_1}=a_1| A_{t-\delta_2}=a_2, Y^0_{t-}, L_{t-},Y^0_{t+s})\\
&&\qquad\rightarrow P(A_{t}=a_1| A_{t-\delta_2}=a_2, Y^0_{t-},
L_{t-},Y^0_{t+s})
\end{eqnarray*}
and
\[
P(A_{t+\delta_1}=a_1| A_{t-\delta_2}=a_2, Y^0_{t-}, L_{t-}) \rightarrow P(A_{t}=a_1| A_{t-\delta_2}=a_2, Y^0_{t-}, L_{t-}),
\]
uniformly in $\delta_2$, as $\delta_1 \downarrow 0$. When the
denominators are bounded from below by a positive number, the ratio
also converges uniformly.

If $a_1 \neq a_2$, then, by the fourth regularity condition, we have
\begin{eqnarray*}
&&\frac{P(A_{t+\delta_1}=a_1| A_{t-\delta_2}=a_2, Y^0_{t-}, L_{t-},Y^0_{t+s})}{\delta_1+\delta_2} \\
&&\qquad\rightarrow \frac{P(A_{t}=a_1| A_{t-\delta_2}=a_2, Y^0_{t-},
L_{t-},Y^0_{t+s})}{\delta_2}
\end{eqnarray*}
and
\begin{eqnarray*}
&&\frac{P(A_{t+\delta_1}=a_1| A_{t-\delta_2}=a_2, Y^0_{t-}, L_{t-})}
{\delta_1+\delta_2}\\
&&\qquad \rightarrow \frac{P(A_{t}=a_1| A_{t-\delta_2}=a_2, Y^0_{t-}, L_{t-})}{\delta_2},
\end{eqnarray*}
uniformly in $\delta_2$, as $\delta_1 \downarrow 0$.\vspace*{1.5pt} Also, the
denominator $\frac{P(A_{t+\delta_1}=a_1| A_{t-\delta_2}=a_2, Y^0_{t-},
L_{t-})}{\delta_1+\delta_2}$ is bounded from below by a positive
number. Hence, we establish the uniform convergence of the ratio.

Combining the two cases above, $|g(\delta_1, \delta_2) -
g_1(\delta_2)|$ is bounded by $O(\delta_1),$ which does not depend on
$\delta_2$, so $g(\delta_1, \delta_2) \rightarrow g_1(\delta_2)$
uniformly in $\delta_2$. Therefore,
\[
\lim_{\delta_2 \downarrow 0} \lim_{\delta_1 \downarrow 0} g(\delta_1, \delta_2)
= \lim_{\delta_1 \downarrow 0} \lim_{\delta_2 \downarrow 0} g(\delta_1, \delta_2).
\]
By the argument at the beginning of the proof, we have proven the first
part of the theorem.

To show that (\ref{eq:markovcon}) implies FTSR, without of loss of
generality, we consider
\begin{eqnarray*}
&&\hspace*{-4.6pt}P\bigl(A_t|L_{t-}, Y^0_{t-}, A_{t-m}, L_{(t-m)-}, Y^0_{(t-m)-}, Y^0_{t+s}\bigr) \\
&&\qquad\hspace*{-4.6pt}=\frac{f(A_t, L_{t-}, Y^0_{t-}, A_{t-m}, L_{(t-m)-}, Y^0_{(t-m)-},
Y^0_{t+s})}
{\sum_{i=0,1}f(A_t=i,L_{t-}, Y^0_{t-}, A_{t-m}, L_{(t-m)-}, Y^0_{(t-m)-}, Y^0_{t+s})}\\
&&\qquad\hspace*{-4.6pt}=\bigl(f(Y^0_{t+s}|A_t, L_{t-}, Y^0_{t-})f\bigl(A_t, L_{t-}, Y^0_{t-},
A_{t-m}, L_{(t-m)-}, Y^0_{(t-m)-}\bigr)\bigr)\\
&&\qquad\quad\hspace*{-4.6pt}
{}\big/\biggl(\sum_{i=0,1}f(Y^0_{t+s}|A_t=i, L_{t-}, Y^0_{t-})\\
&&\qquad\hspace*{39pt}{}\times f\bigl(A_t=i, L_{t-}, Y^0_{t-}, A_{t-m}, L_{(t-m)-}, Y^0_{(t-m)-}\bigr)\biggr)\\
&&\qquad\hspace*{-4.6pt}=\frac{f(Y^0_{t+s}|L_{t-}, Y^0_{t-})f(A_t, L_{t-}, Y^0_{t-}, A_{t-m},
L_{(t-m)-}, Y^0_{(t-m)-})}
{\sum_{i=0,1}f(Y^0_{t+s}|L_{t-}, Y^0_{t-})f(A_t=i, L_{t-}, Y^0_{t-}, A_{t-m}, L_{(t-m)-}, Y^0_{(t-m)-})}\\
&&\qquad\hspace*{-4.6pt}=\frac{f(A_t, L_{t-}, Y^0_{t-}, A_{t-m}, L_{(t-m)-}, Y^0_{(t-m)-})}
{\sum_{i=0,1}f(A_t=i, L_{t-}, Y^0_{t-}, A_{t-m}, L_{(t-m)-}, Y^0_{(t-m)-})}\\
&&\qquad\hspace*{-4.6pt}=P(A_t| L_{t-}, Y^0_{t-}, A_{t-m}, L_{(t-m)-}, Y^0_{(t-m)-}).
\end{eqnarray*}
The second equality follows because of the Markov property. The third
equality uses equation (\ref{eq:markovcon}). We have thus proven the
second half of the theorem.

\section{Simulation parameters}\label{appendixe}
In all simulation models from M1 to M4, we specify the parameters as
follows:
\begin{longlist}
    \item[$\bullet$] let $g(V,t) = C$, a constant; let $C=100$;
    \item[$\bullet$] for M1 (also for M3 and M4), let $\theta = 0.2$ and $\sigma = 1$;
    \item[$\bullet$] for M2, let $m=2$, $\theta_1=0.2$, $\sigma_1 = 1$ and
      $\theta_2=1$, $\sigma_2=0.5$.
      The transition probability of $J_t$ would be $P(t) = e^{At}$,
      where $A=\left({{-1\enskip 1}\atop {1\enskip -1}}\right)$;
    \item[$\bullet$] for initial value, $e_0$ is generated from $N(0,\frac{\sigma}{\sqrt{2\theta}})$;
    \item[$\bullet$] the causal parameter $\Psi = 1$;
    \item[$\bullet$] in M1, M2 and M3, $s(A_t,Y_t) = e^{\alpha_0 + \alpha_1 A_t +
    \alpha_2 Y_t + \alpha_3 A_t Y_t }$; let $\alpha_1 = -0.3$, $\alpha_2
    = -0.005$, $\alpha_3 = 0.007$ and $\alpha_0 = -0.2$;
    \item[$\bullet$] in M4, $A_t$ is generated as follows:
      if $Y_{t-0.5} > 101$ and $Y_{t} > 101$, $s(A_t=1, L^*_t) = 2.8$; if $Y_{t_0.5} < 99$
      and $Y_{t} < 99$, $s(A_t=0, L^*_t) = 2.8$; otherwise, $A_t$ is generated following a model  similar to that in M1, except that $s(A_t,L^*_t)
      =e^{\alpha_0 + \alpha_1 A_t + \alpha_2 L^*_t + \alpha_3 A_t
      L^*_t}$; the values of the $\alpha$'s are the same as before;
    \item[$\bullet$] in M4, $\eta_t$ follows an Ornstein--Uhlenbeck process with
    parameters $\theta = 0.2$ and $\sigma = 1$;
    \item[$\bullet$] for initial value, $A_0$ is generated from $\operatorname{Bernoulli}(\operatorname{expit}(\alpha_0+\alpha_2 Y_0))$;
    \item[$\bullet$] $K=5$ is the number of periods;
    \item[$\bullet$] number of subjects $n=5000$.
\end{longlist}
\end{appendix}

\section*{Acknowledgments} The authors would like to thank the
causal inference reading group at the University of Pennsylvania and Judith Lok for helpful
discussions.

\printaddresses

\end{document}